\documentclass[11pt, reqno]{amsart}
\usepackage{mathtools,bm}
\usepackage{amsmath}
\usepackage{amsthm}
\usepackage{amssymb}
\usepackage{mathrsfs}
\usepackage{verbatim}
\usepackage{esint}
\usepackage{pdfpages}
\usepackage[foot]{amsaddr}
\usepackage[english]{babel}

\usepackage[margin=1.00in]{geometry}

\usepackage{esint}

\usepackage[all]{xy}
\catcode`\@=11
\def\theequation{\@arabic{\c@section}.\@arabic{\c@equation}}
\catcode`\@=12

\usepackage{subfigure}
\usepackage{xcolor}
{%
	\setcounter{enumi}{0}

	\begin{enumerate}}%
	{\end{enumerate} }

{%
	\setcounter{enumi}{0}

	\begin{enumerate}}%
	{\end{enumerate} }

\usepackage[colorlinks]{hyperref}
\definecolor{myblue}{rgb}{0.0, 0.0, 1.0}
\definecolor{mygreen}{rgb}{0.2,0.8,0.2}
\hypersetup{
	colorlinks=true,
	linkcolor=myblue,         
	citecolor=mygreen,
	urlcolor=mygreen,
	linktocpage=true
}
\usepackage[nameinlink]{cleveref}
\usepackage[hyperpageref]{backref}
\newtheorem{theorem}{Theorem}[section]
\newtheorem{lemma}[theorem]{Lemma}
\newtheorem{proposition}[theorem]{Proposition}

\theoremstyle{definition}
\newtheorem{definition}[theorem]{Definition}
\newtheorem{remark}[theorem]{Remark}

\numberwithin{equation}{section}

\newcommand*\re{\mathbb{R}}
\newcommand*\rd{\mathbb{R}^d}

\newcommand*\N{\mathcal{N}}

\newcommand{\ep} {\epsilon}
\newcommand{\al} {\alpha}
\newcommand{\pa} {\partial}
\newcommand{\be} {\beta}
\newcommand{\de} {\delta}
\newcommand{\De} {\Delta}

\newcommand{\ga} {\gamma}

\newcommand{\Om} {\Omega}
\newcommand{\la} {\lambda}

\newcommand{\no} {\nonumber}
\newcommand{\noi} {\noindent}

\newcommand{\ra} {\rightarrow}

\def\dx{{\rm d}x}
\def\dy{{\rm d}y}
\def\dt{{\rm d}t}

\def\dtau{{\rm d}\tau}

\def\C{{\mathcal C}}

\def\M{{\mathcal M}}
\def\S{{\mathcal S}}

\def\R{{\mathbb R}}
\def\N{{\mathbb N}}
\def\F{{\mathcal F}}
\def\R{{\mathbb R}}

\def\N{{\mathcal N}}
\def\F{{\mathcal F}}
\def\({{\Big(}}
\def\){{\Big)}}

\def\ws2{{\F_{\frac{N}{2}}}}
\def\l2{{ L^{1,\;\infty}(\log L)^2}}
\def\M2{\mathcal M_{\log L}}
\def\cc{{\C_c^\infty}}

\def\p{{p^{\prime}}}

\def\q{{q^{\prime}}}
\def\dt{{\rm d}t}
\def\dx{{\rm d}x}
\def\dz{{\rm d}z}

\def\dmuo{{\rm d}\mu_1}
\def\dmut{{\rm d}\mu_2}

\def\wrs{{W_0^{s,r}(\Om)}}
\def\wpso{{W_0^{s_1,p}(\Om)}}
\def\wqst{{W_0^{s_2,q}(\Om)}}

\newcommand{\intom}{\displaystyle{\int_{\Omega}}}

\DeclarePairedDelimiter\abs{\lvert}{\rvert}%
\DeclarePairedDelimiter\norm{\lVert}{\rVert}%
\makeatletter
\let\oldnorm\norm
\def\norm{\@ifstar{\oldnorm}{\oldnorm*}}

\DeclareMathOperator*{\lowlim}{\underline{lim}}

\newcommand{\normo}{[u]_{s_1,p}}
\newcommand{\normt}{[u]_{s_2,q}}

\newcommand{\Tail} {\mathrm{Tail}}
\def \tb {\textcolor{blue}}

\newcommand{\II}{\iint\limits_}
\newcommand{\lao}{\lambda_{s_1,p}^1}
\newcommand{\lat}{\lambda_{s_2,q}^1}

\allowdisplaybreaks

\title[On generalized nonlocal eigenvalue problems]{On generalized eigenvalue problems of fractional $(p,q)$-Laplace operator with two parameters}

\author[Biswas and Sk]{Nirjan Biswas$^*$ and Firoj Sk}

\address{N. Biswas, Tata Institute of Fundamental Research, Centre For Applicable Mathematics, Post Bag No 6503, Sharada Nagar, Bangalore 560065, India.}

\address{F. Sk, Carl von Ossietzky Universit\"{a}t Oldenburg,
Fakult\"{a}t V, Institut f\"{u}r Mathematik, Ammerländer Heerstraße 114-118, 26129 Oldenburg, Germany.}
\email{nirjan22@tifrbng.res.in, firoj.sk@uol.de}
\thanks{$^*$Corresponding author}
\smallskip


\begin{document}

\maketitle

\begin{abstract}
    For $s_1,s_2\in(0,1)$ and $p,q \in (1, \infty)$, we study the following nonlinear Dirichlet eigenvalue problem with parameters $\alpha, \beta \in \mathbb{R}$ driven by the sum of two nonlocal operators:
    \begin{equation*}
        (-\Delta)^{s_1}_p u+(-\Delta)^{s_2}_q u=\alpha|u|^{p-2}u+\beta|u|^{q-2}u\;\;\text{in }\Omega, \quad u=0\;\;\text{in } \mathbb{R}^d \setminus \Omega, \ \ \ \qquad \quad \mathrm{(P)}
    \end{equation*}
where $\Omega \subset \mathbb{R}^d$ is a bounded Lipschitz open set. Depending on the values of $\alpha,\beta$, we completely describe the existence and non-existence of positive solutions to (P). We construct a continuous threshold curve in the two-dimensional $(\alpha, \beta)$-plane, which separates the regions of the existence and non-existence of positive solutions. In addition, we prove that the first Dirichlet eigenfunctions of the fractional $p$-Laplace and fractional $q$-Laplace operators are linearly independent, which plays an essential role in the formation of the curve. Furthermore, we establish that every nonnegative solution of (P) is globally bounded.
\end{abstract}
\tableofcontents

\section{Introduction and Main Results}

In this paper, we are concerned with the existence and non-existence of positive solutions to the following nonlinear eigenvalue problem involving the fractional $(p,q)$-Laplace operator with zero Dirichlet boundary condition:   

\begin{equation}
\tag{\text{EV};\;$\alpha,\beta$}\label{main problem}
    (-\Delta)^{s_1}_p u + (-\Delta)^{s_2}_q u = \alpha|u|^{p-2}u + \beta|u|^{q-2}u 
    \;\;\text{in }\Omega, \quad 
    u=0 \;\;\text{in } \rd \setminus \Om,
\end{equation}
where $0 < s_2 < s_1 < 1 < q < p < \infty$, $\al, \be \in \re$ are two parameters and $\Omega\subset\rd$ is a bounded Lipschitz open set. In general, the fractional $r$-Laplacian $(-\De)^{s}_r$ ($s \in (0,1)$ and $r \in (1, \infty)$) is defined as
$$
(-\Delta)^{s}_r u(x) := \text{P.V.}\int_{\rd}\frac{|u(x)-u(y)|^{r-2}(u(x)-u(y))}{|x-y|^{d+sr}} \,\dy, \;\;\;x\in \rd,
$$
where P.V. stands for the principle value. 

The local counterpart of \eqref{main problem} is the following Dirichlet eigenvalue problem for the $(p,q)$-Laplace operator:
\begin{equation}\label{local EV}
    -\Delta_p u - \Delta_q u = \alpha|u|^{p-2}u+\beta|u|^{q-2}u\;\;\text{in }\Omega, \quad 
    u=0\;\;\text{in } \partial \Om.
\end{equation}
The study of $(p,q)$-Laplace operators are well known for their applications in physics, chemical reactions, reaction-diffusion equations e.t.c. for details, see \cite{ChIl, De,  Fi} and the references therein. Some authors considered the eigenvalue problems for the $(p,q)$-Laplace operator. In this direction, for $\al=\be$, Motreanu-Tanaka in \cite{MoTa} obtained the existence and non-existence of positive solutions of \eqref{local EV}. For $\al \neq \be$, in \cite{BoTa1} Bobkov-Tanaka extended this result by providing a certain region in the $(\al,\be)$-plane that allocates the sets of existence and non-existence of positive solutions of \eqref{local EV}. Moreover, they constructed a \textit{threshold curve} in the first quadrant of the $(\al,\be)$-plane, which separates these two sets. Later, in \cite{BoTa2}, the same authors plotted a different curve for the existence of ground states and the multiplicity of the positive solutions for \eqref{local EV}. It is essential that in which region the positive solution of \eqref{local EV} exists or does not exist, and the behaviour of the threshold curve depends on whether $\phi_{p}, \phi_{q}$ are linearly independent, where $\phi_{p} \text{ and } \phi_{q}$ are the first Dirichlet eigenfunctions of the operators $-\Delta_p$ and $-\Delta_q$ respectively.  For other results related to the positive solutions of eigenvalue problems involving $(p,q)$-Laplace operator, we refer to \cite{BaPaZe, BoTa3, Ta1} and the references therein.  

In the nonlocal case, parallelly, many authors studied the nonlinear equations driven by the sum of fractional $p$-Laplace and fractional $q$-Laplace operators with the critical exponent. For example, see \cite{Am1, AmIs, BhMu, GiKuSr, GoKuSr} where the weak solution's existence, regularity, multiplicity, positivity and other qualitative properties are investigated. The study of \eqref{main problem} is motivated by the Dancer-Fu\v{c}ik (DF) spectrum of the fractional $r$-Laplace operator. The DF spectrum  of the operator $(-\De)^{s}_r$ is the set of all points $(\al,\beta)\in\re^2$ such that the following problem
\begin{equation}\label{Fucik}
(-\Delta)_r^s u = \al(u^+)^{r-1}-\be(u^-)^{r-1} \text{ in }\Om, \quad u = 0\text{ in } \rd \setminus \Om,
\end{equation}
admits a nontrivial weak solution, where $u^{\pm} = \max\{\pm u, 0 \}$ is the positive and negative part of $u$. For $r=2$, in \cite{GoSr}, Goyal-Sreenadh considered \eqref{Fucik} and proved the existence of a first nontrivial curve in the DF spectrum. They also showed that the curve is Lipschitz continuous, strictly decreasing, and studied its asymptotic behaviour. For $r \neq 2$, in \cite{PeSqYa}, the authors constructed an unbounded sequence of decreasing curves in the DF spectrum. Nevertheless, the study of the spectrum for the fractional $(p,q)$-Laplace operator is not well explored. In \cite{NgVo}, for $\al=\be$, Nguyen-Vo studied the following weighted eigenvalue problem with zero Dirichlet boundary condition:
\begin{equation}\label{weight problem1}
   (-\Delta)^{s_1}_pu + (- \Delta)^{s_2}_qu = \alpha\left(m_p |u|^{p-2}u + m_q|u|^{q-2}u \right) \;\;\text{in }\Omega, \quad u = 0\;\;\text{in } \rd \setminus \Om,
\end{equation}
where $0< s_2< s_1< 1< q \le p< \infty$, the weights $m_p, m_q$ are bounded in $\Omega$ and satisfy $m_p^+, m_q^+ \not \equiv 0$. Depending on the values of $\al$, the authors obtained the existence and non-existence of positive solutions of \eqref{weight problem1}. 

The primary aim of this paper can be summarized into the following two aspects:

\noi \textbf{(a)} We provide a comprehensive analysis of the sets in the $(\alpha, \beta)$-plane that determine the existence and non-existence of positive solutions for the equation \eqref{main problem}. Following the local case approach, we construct a continuous threshold curve denoted as $\mathcal{C}$ that effectively separates the regions where positive solutions exist from those where they do not. In some specific regions of the $(\alpha, \beta)$-plane, we employ the sub-super solutions technique to establish the existence of positive solutions. To apply this technique, we utilize the crucial result stated in Theorem \ref{bounded solution}, which proves that every nonnegative solution of \eqref{main problem} is globally bounded.

\noi \textbf{(b)} The existence and non-existence of positive solutions to \eqref{main problem} depend on the following statement:
\begin{equation}\label{LI}
  \tag{LI}  \phi_{s_1,p} \neq c \phi_{s_2,q} \text{ for any } c \in \R,
\end{equation}
where $\phi_{s_1,p} \text{ and } \phi_{s_2,q}$ are the first eigenfunctions of the operators $(-\Delta)_p^{s_1}$ and $(-\Delta)_q^{s_2}$ corresponding to the first eigenvalues $\lao$ and $\lat$ respectively in $\Omega$ under zero Dirichlet boundary condition. While this linear independence condition for the operators $-\Delta_p$ and $-\Delta_q$ was conjectured in \cite{BoTa1} and later proved in \cite{BoTa2}, its validity remains unknown for any $s_1, s_2 \in (0,1)$. Nevertheless, several authors have assumed the condition \eqref{LI} in various contexts (e.g., \cite{GiGoMo, NgVo}). We establish the validity of \eqref{LI} under certain assumptions on $s_1$ and $s_2$, as demonstrated in Theorem \ref{linearindependent}.

Recall that, for $0<s<1\leq r<\infty$, the fractional Sobolev space is defined as
$$W^{s,r}(\Omega):=\left\{u\in L^r(\Omega):[u]_{s,r,\Omega}<\infty\right\},$$
with the so-called fractional Sobolev norm $\norm{u}_{s,r,\Omega} :=(\|u\|_{L^r(\Omega)}^r + [u]_{s,r,\Omega}^r)^\frac{1}{r}, $
where
$$ [u]_{s,r,\Omega}^r :=\iint\limits_{\Omega\times\Omega}\frac{|u(x)-u(y)|^r}{|x-y|^{d+sr}}\,\dx\dy, $$
is called the Gagliardo seminorm. For $r \in (1,\infty)$, $W^{s,r}(\Omega)$ is a reflexive Banach space with respect to the fractional Sobolev norm $\norm{\cdot}_{s,r,\Om}$. 
Now we consider the following closed subspace of $W^{s,r}(\rd)$:
$$W_0^{s,r}(\Omega) :=\{u\in W^{s,r}(\rd) : u=0 \text{ in }\rd \setminus \Om\},$$
endowed with the seminorm $[\cdot]_{s,r,\rd}$, which is an equivalent norm in $W_0^{s,r}(\Omega)$ (\cite[Lemma 2.4]{BrLiPa}). For details of the fractional Sobolev spaces and their related embedding results, we refer to \cite{BrLiPa, BrPa, DiPaVa} and the references therein. For $s_1 > s_2$ and $p > q \geq 1$, the continuous embedding $\wpso \hookrightarrow \wqst$ (see \cite[Proposition~2.2]{AnWa}) allows us to introduce the notion of weak solution for \eqref{main problem} in the following sense:

\begin{definition}\label{weak formulation of frac p q Lap}
A function $u\in\wpso$ is called a weak solution of \eqref{main problem} if the following identity holds for all $\phi \in \wpso$:
\begin{multline*}
\iint\limits_{\rd\times\rd}\frac{|u(x)-u(y)|^{p-2}(u(x)-u(y))(\phi(x)-\phi(y))}{|x-y|^{d+s_1p}} \, \dx \dy
     \\
    + \iint\limits_{\rd\times\rd}\frac{|u(x)-u(y)|^{q-2}(u(x)-u(y))(\phi(x)-\phi(y))}{|x-y|^{d+s_2q}} \, \dx \dy
    = \al\intom |u|^{p-2}u\phi\,\dx+\be\intom|u|^{q-2}u\phi\,\dx.
\end{multline*}
\end{definition}

In our first theorem, we prove the existence of a positive solution for \eqref{main problem} if any of $\al$ and $\be$ is larger than the first Dirichlet eigenvalue of the fractional $p$-Laplacian and fractional $q$-Laplacian respectively. We also show that this range of $\al, \be$ is necessary for the existence of a positive solution when \eqref{LI} does not hold.

\begin{figure}[htp]\label{fig}
    \centering
    \subfigure[The case (LI) holds (with $\al^*_{s_1,p} < \infty$)]{
    \includegraphics[width=9cm]{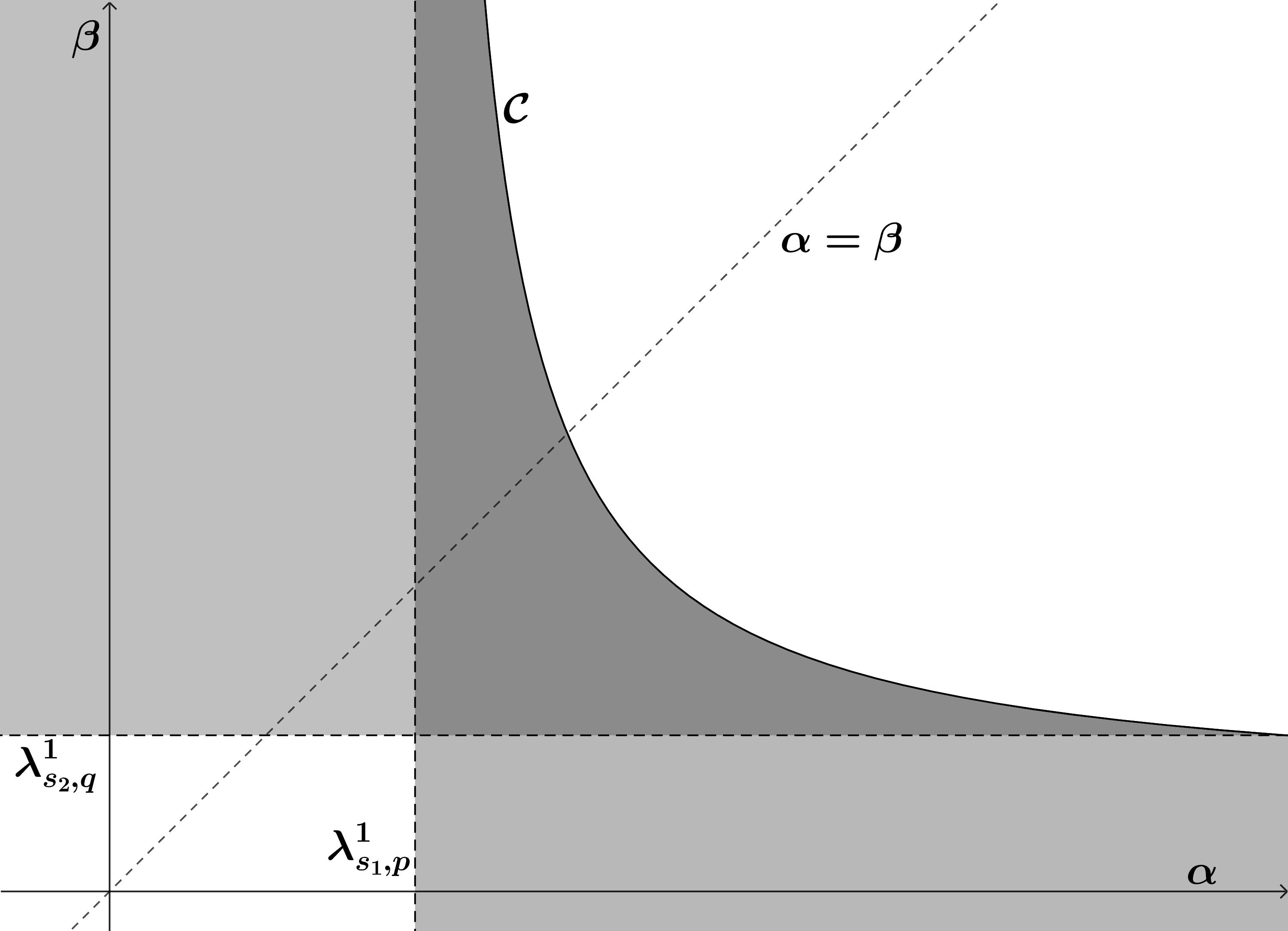}
    }
    \subfigure[The case (LI) does not hold]{
    \includegraphics[width=9cm]{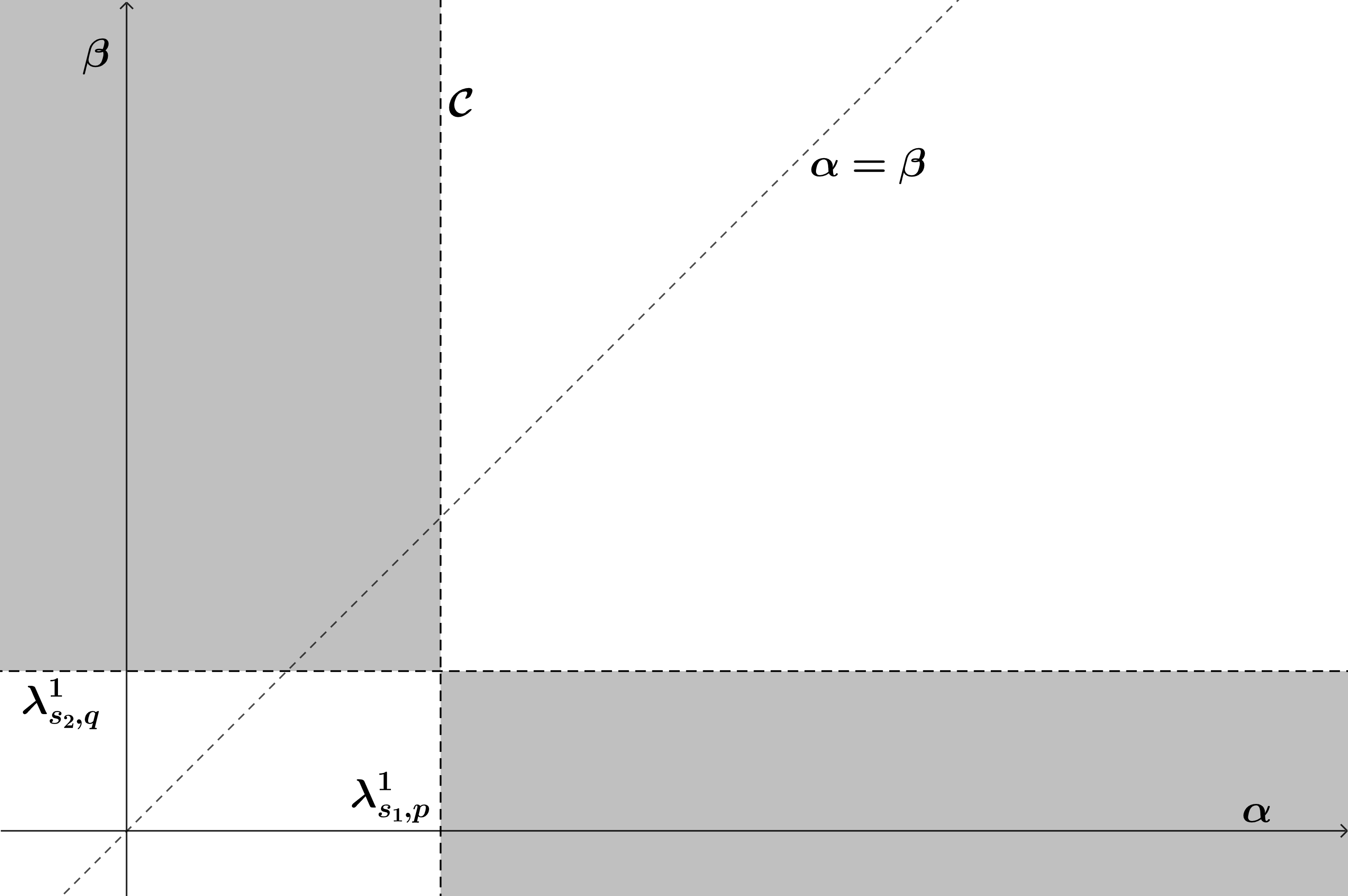} 
    }
    \caption{Shaded region denotes existence, and unshaded region denotes non-existence of positive solutions.}
\end{figure}

\begin{theorem}\label{Existence and nonexistence 1} 
Let $0 < s_2 < s_1 < 1 < q < p < \infty$. Assume that 
\begin{equation}\label{assumption1}
(\al, \be) \in \left((\lao, \infty) \times (-\infty, \lat)\right) \cup \left((-\infty, \lao) \times  (\lat,\infty)\right) \cup \left( \{\lao\} \times \{\lat\} \right).
\end{equation}
The following hold (see Figure \ref{fig}):
\begin{itemize}
    \item[(i)] \rm{{(\textbf{Sufficient condition})}}: Let $\al, \be$ satisfy  \eqref{assumption1}. In the case, when $\al= \lao$ and $\be = \lat$, we assume that \eqref{LI} violates. Then \eqref{main problem} admits a positive solution.
    \smallskip
    \item[(ii)] \rm{{(\textbf{Necessary condition})}}: Let \eqref{LI} violates and \eqref{main problem} admits a positive solution. Then $\al, \be$ satisfy  \eqref{assumption1}.
\end{itemize}
\end{theorem}

\begin{remark}
(i) The above theorem asserts that (\tb{EV; $\lao,\,\lat$})  admits a positive solution if and only if \eqref{LI} violates. Indeed, (\tb{EV; $\lao,\,\lat$}) admits a non-trivial solution only when \eqref{LI} violates (see Proposition \ref{Nonexistence1}). 

\noi (ii) If \eqref{LI} violates, then Theorem \ref{Existence and nonexistence 1} gives a complete description of the set of existence and non-existence of positive solutions of \eqref{main problem}. In particular, Theorem \ref{Existence and nonexistence 1} generalizes the result of \cite[Theorem 1.1]{NgVo} for $\al \neq \be$. 
\end{remark}

It is observed that for $\al,\be\in\re$, the problem \eqref{main problem} is equivalent to the problem (\tb{EV; $\be + \theta, \be$}), where $\theta = \al-\be.$ Using this terminology we define the following curve:
\begin{definition}[Threshold curve]
For brevity, denote $\be = \la$. For each $\theta \in \re$ consider the following quantity:
\begin{equation}\label{lambda* definition}
    \la^*(\theta) := \sup\left\{\la \in \re : (\tb{\text{EV; $\la+\theta,\la$}})\text{ has a positive solution}\right\}.
\end{equation}
If such $\la$ does not exist, we then set $\la^*(\theta) = -\infty.$ The \textit{threshold curve} corresponding to \eqref{main problem} is defined as $\mathcal{C} := \{(\la^*(\theta) + \theta, \la^*(\theta)): \theta \in \R\}$. Also, we define the following quantities:
\begin{equation*}\label{theta* definition}
    \theta^* :=\la^1_{s_1,p}-\la^1_{s_2,q}, \, \al^*_{s_1,p} := \frac{[\phi_{s_2,q}]_{s_1,p, \rd}^p}{\norm{\phi_{s_2,q}}_{L^p(\Om)}^p}, \text{ and }\theta^*_{+} := \al^*_{s_1,p} -\la^1_{s_2,q}. 
\end{equation*}
Clearly, $\theta^* \leq \theta^*_+$ and $\theta^* = \theta^*_+$ if and only if \eqref{LI} violates (from (iv) of Proposition \ref{first eigenvalue}). 
\end{definition}

In the following proposition, we discuss some qualitative properties of $\mathcal{C}$ and see that $\mathcal{C}$ carries similar behaviours as in the local case \cite[Proposition 3 and Fig. 2]{BoTa1}.

\begin{proposition}\label{threshold curve}
Let $0 < s_2 < s_1 < 1 < q < p < \infty$. Then the following hold:
\begin{itemize}
    \item [(i)]$\la^*(\theta) < \infty$ for all $\theta\in\re$.
    \smallskip
    \item[(ii)] $\la^*(\theta^*) + \theta^* > \la^1_{s_1,p}$ and $\la^*(\theta^*) > \la^1_{s_2,q}$ if and only if \eqref{LI} holds.
    \smallskip
    \item[(iii)] $\la^*(\theta) + \theta \geq \la^1_{s_1,p}$ and $ \la^*(\theta) \geq \la^1_{s_2,q}$ for all $\theta \in \re.$
    \smallskip
    \item[(iv)]$\la^*(\theta)$ is decreasing and $\la^*(\theta) + \theta$ is increasing on $\re.$
    \smallskip
    \item[(v)] If $\al^*_{s_1,p}$ is finite, then $\la^*(\theta) = \la^1_{s_2,q}$ for all $\theta \geq \theta^*_+.$
    \smallskip
    \item[(vi)]$\la^*$ is continuous on $\re.$
\end{itemize}
\end{proposition}

According to (iii) of the above proposition, $\mathcal{C} \subset ([\lao, \infty) \times [\lat, \infty))$. Further, if $\al^*_{s_1,p} = \infty$, from the property (iii), we observe that $\mathcal{C}$ always lies above the line $\be = \lat$. From now onwards, we assume that $\al^*_{s_1, p} < \infty$. 
In the following theorem, we demonstrate that $\mathcal{C}$ separates the sets of existence and non-existence of positive solutions in the region  $([\lao, \infty) \times [\lat, \infty))$ (see Figure \ref{fig}). 

\begin{theorem}\label{Existence main 1}
Let $0 < s_2 < s_1 < 1 < q < p < \infty$. Let $\al \ge \lao$ and $\be \ge \lat$. Assume that \eqref{LI} holds. 
\begin{itemize}
    \item[(i)] If $\be \in (\lat, \la^*(\theta))$, then \eqref{main problem} admits a positive solution.
    \smallskip
    \item[(ii)] If $\al > \lao$ and $\be < \la^*(\theta)$, then \eqref{main problem} admits a positive solution.
    \smallskip
    \item[(iii)]  If $\be > \la^*(\theta)$, then there does not exist any positive solution for \eqref{main problem}. 
\end{itemize}
\end{theorem}

Now we state the existence and non-existence of positive solutions on the curve $\mathcal{C}$ (see Figure \ref{fig}). 

\begin{theorem}\label{Existence main 2}
Let $0 < s_2 < s_1 < 1 < q < p < \infty$. 
\begin{enumerate}
    \item[(i)] If $\theta < \theta^*_+$, then {\rm (\tb{EV; $\la^*(\theta) + \theta, \la^*(\theta)$})} admits a positive solution.
    \smallskip 
    \item[(ii)] If $\theta > \theta^*_+$, then there does not exist any positive solution for { \rm(\tb{EV; $\la^*(\theta) + \theta, \la^*(\theta)$})}. 
\end{enumerate}
\end{theorem}

The above theorem does not consider the borderline case $\theta = \theta_+^*$. In this case, we have a partial result in Remark \ref{borderline 2}, which says that (\tb{EV; $\la^*(\theta)+\theta,\la^*(\theta)$}) does not admit any ground state solution. 
 
\begin{remark}\label{range changes}
    The relations among $s_1, s_2, p, q$ are taken without loss of any generality. All the preceding results in this paper hold for the remaining cases by choosing the appropriate solution space as given below:
\begin{enumerate}
    \item[(i)]  For $s_1 < s_2$ and $p < q$ (symmetric), we choose the solution space as $\wqst$.
    \item[(ii)] For $s_2 < s_1$  and $p<q$ (cross), we choose 
   the solution space as $\wpso \cap \wqst$ endowed with the norm $[\cdot]_{s_1,p,\rd}+[\cdot]_{s_2,q,\rd}$.
   \item[(iii)] For $s_1 = s_2=s$ and $p \neq q$, we choose the solution space as $W^{s,p}_0(\Omega) \cap W^{s,q}_0(\Omega)$ endowed with the norm $[\cdot]_{s,p,\rd}+[\cdot]_{s,q,\rd}$.
   \end{enumerate}
\end{remark}

The next theorem verifies the linear independency of the first Dirichlet eigenfunctions of the fractional $p$-Laplacian and the fractional $q$-Laplacian.

\begin{theorem}\label{linearindependent}
Let $1< q< p<\infty$ and $s_1, s_2 \in (0,1)$ satisfy the following condition:
\begin{align*}
    \frac{s_1\p}{\q}<s_2<s_1.
\end{align*}
Then the set $\{\phi_{s_1, p} , \phi_{s_2, q}\}$ is linearly independent. 
\end{theorem}

\begin{remark}
    Theorem \ref{linearindependent} holds if we take the other relations among $s_1, s_2, p, q$ listed below:
    \begin{enumerate}
        \item[(i)] For  $1< p<q<\infty$ and $\frac{s_2 \q}{\p} < s_1 <s_2$ (interchanging the roles of $s_1, s_2, p, q$). 
        \item[(ii)] For $1 < q<p<\infty$ and $s_1=s_2$. 
    \end{enumerate}
\end{remark}

The rest of the paper is organized as follows. \Cref{section2} briefly discusses the first Dirichlet eigenpair of fractional $r$-Laplace operator, recalls the discrete Picone's inequalities, and proves some technical results. In \Cref{section5}, we prove the validity of \eqref{LI}. This section contains the proof of Theorem \ref{linearindependent}. In \Cref{section6}, we establish the regularity of the solution for \eqref{main problem} and state a version of the strong maximum principle related to 
\eqref{main problem}.  \Cref{section3} studies various frameworks of energy functionals associated with \eqref{main problem}. Finally, \Cref{Section4} studies the existence and non-existence of positive solutions for \eqref{main problem}. In this section, we prove Theorem \ref{Existence and nonexistence 1}-\ref{Existence main 2} and Proposition \ref{threshold curve}. 

\section{Preliminaries}\label{section2}
In this section, we recall some qualitative properties of the first nonlocal eigenvalue and its corresponding eigenfunction. Afterwards, we recall the discrete Picone's identities. We 
list the following notations to be used in this paper:

\noi \textbf{Notation:} 
\begin{itemize}
    \item $B_R(x)\subset\rd$ denotes an open ball of radius $R>0$ centered at $x.$
    \smallskip
    \item For a set $E\subset\rd$, $|E|$ denotes the Lebesgue measure of $E$.
    \smallskip
    \item We denote $\dmuo:= \abs{x-y}^{-(d+s_1p)} \dx \dy $ and $\dmut:= \abs{x-y}^{-(d+s_2q)} \dx \dy $.
    \smallskip
    \item For $r \in (1, \infty)$, the conjugate of $r$ is denoted as $r^{\prime}:=\frac{r}{r-1}$.
    \smallskip
    \item For $0 < s < 1 < r < \infty$, we denote $[\cdot]_{s,r,\rd}$ as $[\cdot]_{s,r}$, and $\norm{\cdot}_{L^r(\Om)}$ as $\norm{\cdot}_r$.
    \smallskip
    \item For $k \in \mathbb{N}$, we denote $u_k(x) := u(x) + \frac{1}{k}$ where $x \in \rd$.
    \smallskip
    \item For $\ga \in (0,1)$, the H\"{o}lder seminorm $[f]_{\C^{0, \ga}(\Om)} := \underset{x,y \in \Om,\, x \neq y}{\sup} \displaystyle \frac{\abs{f(x)- f(y)}}{\abs{x-y}^{\ga}}$.
    \smallskip
    \item For $sr < d$ (where $0 < s < 1 < r < \infty$), the fractional critical exponent $\displaystyle r^*_s := \frac{rd}{d-sr}$.
    \smallskip
    \item For each $n \in \mathbb{N}$, we denote the positive and negative parts $(f_n)^{\pm}$ by $f_n^{\pm} := \max\{\pm f_n, 0\}$.
    \smallskip
    \item Eigenvalue of \eqref{nonlocal Dirichlet EV}, $\lambda_{s,r}(\Om)$ is denoted as $\lambda_{s,r}$.
    \smallskip
    \item We denote the eigenfunction of \eqref{nonlocal Dirichlet EV} corresponding to the first eigenvalue $\la^1_{s, r}$ as $\phi_{s,r}$.
    \smallskip
    \item For $r \in (1, \infty)$, $x_0 \in \Om$ and $R>0$, the nonlocal tail of $f \in W^{s,r}_0(\Om)$ is defined as \begin{align*}
        \Tail_r(f;x_0,R) := \left(R^{sr}\int_{\rd\setminus B_R(x_0)}\frac{|f(x)|^{r-1}}{|x-x_0|^{d+sr}}\,\dx\right)^{\frac{1}{r-1}}.
    \end{align*}
    \smallskip
    \item $C$ is denoted as a generic positive constant.
\end{itemize}

\subsection{First eigenvalue of fractional $r$-Laplacian}
For a bounded open set $\Om \subset \rd$ and $0 < s < 1 < r < \infty$, we consider the following nonlinear eigenvalue problem: 
\begin{equation}\label{nonlocal Dirichlet EV}
(-\Delta)_r^su = \lambda_{s,r}|u|^{r-2}u\;\text{in }\Omega, \quad 
u = 0\;\text{ in } \rd \setminus \Om.
\end{equation}
We say $\la_{s,r}$ is an eigenvalue of \eqref{nonlocal Dirichlet EV}, if there exists non-zero $u \in \wrs$ satisfying the following identity for all $\phi\in \wrs$:
\begin{equation*}
    \begin{split}
    \iint\limits_{\rd\times\rd} \frac{|u(x)-u(y)|^{r-2}(u(x)-u(y))(\phi(x)-\phi(y))}{|x-y|^{d+sr}} \,  \dx \dy = \la_{s,r} \intom |u(x)|^{r-2} u(x) \phi(x) \, \dx.
    \end{split}
\end{equation*}
In this case, $u$ is called an eigenfunction corresponding to  $\lambda_{s,r}$, and we denote $\left( \lambda_{s,r}, u \right)$ as an eigenpair. In the following proposition, we collect some qualitative properties of the first eigenpair of \eqref{nonlocal Dirichlet EV}.
\begin{proposition}\label{first eigenvalue}
For $r \in (1, \infty)$ and $s \in (0,1)$, consider the following quantity:
\begin{align*}
    \lambda^1_{s,r} = \inf \left\{ [u]_{s,r}^r : u \in \wrs \text{ and } \intom|u|^r=1 \right\}.
\end{align*}
Then the following hold:
\begin{itemize}
    \item[(i)] $\la^1_{s,r}$ is the first positive eigenvalue of \eqref{nonlocal Dirichlet EV}. 
    \smallskip
    \item[(ii)] Every eigenfunction corresponding to $\la^1_{s,r}$ has a constant sign in $\Om$.
    \smallskip
    \item[(iii)] If $v$ is an eigenfunction of \eqref{nonlocal Dirichlet EV} corresponding to an eigenvalue $\la_{s,r}>0$ such that $v>0$ a.e. in $\Om$, then $\la_{s,r}=\la^1_{s,r}$.
    \smallskip
    \item[(iv)] Any two eigenfunctions corresponding to $\lambda^1_{s,r}$ are constant multiple of each other.
    \smallskip
    \item[(v)] Any eigenfunction of \eqref{nonlocal Dirichlet EV} corresponding to an eigenvalue $\lambda_{s,r}$ is in $L^{\sigma}(\rd)$ for every $\sigma \in [1, \infty]$. Moreover, if $\Om$ is of class $\mathcal{C}^{1,1}$, then the eigenfunction lies in $\C^{0, \ga}(\overline{\Om})$ for some $\ga \in (0,s]$. 
    \smallskip
\end{itemize}
\end{proposition}

\begin{proof}
For proof of (i) and (iii), we refer to~\cite[Lemma 2.1 and Theorem 4.1]{FrPa}. For (ii), see~\cite[Proposition 2.6]{BrPa}. Then the proof of (iv) follows using~\cite[Theorem 4.2]{FrPa}. 

(v) Let $u$ be an eigenfunction of \eqref{nonlocal Dirichlet EV} corresponding to $\la_{s,r}$. By~\cite[Theorem 3.3]{BrLiPa}, $u \in L^{\infty}( \rd )$. Further, since $u\in W^{s,r}(\rd) \cap L^\infty(\rd)$, the interpolation argument yields $u \in L^{\sigma}(\rd)$ for every $\sigma \geq r$. Also for $\sigma \in [1, r)$, applying H\" {o}lder's inequality with the conjugate pair $(\frac{r}{\sigma}, \frac{r - \sigma}{\sigma})$,  
\begin{align*}
\int_{\Om} \abs{u}^{\sigma} \le   \left( \int_{\Om} \abs{u}^r \right)^{\frac{\sigma}{r}} \abs{\Om}^{\frac{r - \sigma }{r}}.  
\end{align*}
Thus, $u \in L^{\sigma}(\rd)$ for every $\sigma \in [1, \infty]$.
Furthermore, we apply~\cite[Theorem 1.1]{IaMoSq} to get $u \in C^{0, \ga}(\overline{\Om})$ for some $\ga \in (0,s]$.
\end{proof}

\subsection{Some important results}
In this subsection, we state some elementary inequalities, recall Picone's inequalities for nonlocal operators and collect some test functions in $W^{s,r}_{0}(\Omega)$.

\begin{lemma}\label{inequality}
Let $\ga \in \R^+$. The following hold: 
\begin{enumerate}
    \item[(i)] If $\ga>1$, then 
    \begin{align*}
        &\abs{f(x)-f(y)}^{\ga-2}(f(x)-f(y))(f^+(x) - f^+(y)) \ge \abs{f^+(x) - f^+(y)}^{\ga}; \\
        &\abs{f(x) - f(y)}^{\ga-2}(f(x)-f(y))(f^-(y) - f^-(x)) \ge \abs{f^-(x) - f^-(y)}^{\ga},
    \end{align*}
    where $f^{\pm} = \max \{\pm f,0\}$.
    \item[(ii)] Let $a, b \in \R$. If $\ga \ge 2$, then $\abs{a-b}^{\ga-2}(a-b) \le C \left( \abs{a}^{\ga-2}a - \abs{b}^{\ga-2}b \right)$ for some $C = C(\ga)>0$.
    \item[(iii)] Let $a, b \in \R$. Then $\abs{\abs{a}^{\ga} - \abs{b}^{\ga}} \le \ga \left( \abs{a}^{\ga-1} + \abs{b}^{\ga-1} \right) \abs{a-b}$.
\end{enumerate}
\end{lemma}

\begin{proof}
Proof of (i) follows from \cite[Lemma A.2]{BrPa}. Proof of (ii) follows from \cite[(2.2) of Page 5]{IaMoSq1}.
Proof of (iii) follows using the fundamental theorem of calculus. 
\end{proof}
We recall several versions of the discrete Picone's inequality that are useful in proving our results.
\begin{lemma}[Discrete Picone's inequality]\label{picone}
Let $r_1 , r_2 \in (1 , \infty)$ with $r_2 \leq r_1$ and let $f, g: \rd \to \re$ be two measurable functions with $f > 0, \; g \ge 0$. Then the following hold: \begin{enumerate}
    \item[(i)] For $x, y \in \rd$,
    $$
        |f(x)-f(y)|^{r_1-2}(f(x)-f(y)) \left(\frac{g(x)^{r_2}}{f(x)^{r_2-1}} - \frac{g(y)^{r_2}}{f(y)^{r_2-1}} \right) \le \abs{g(x)-g(y)}^{r_2} \abs{f(x)-f(y)}^{r_1-r_2}. 
    $$
    \item[(ii)] For $x,y \in \rd$,
    \begin{align*}
        |f(x)-f(y)|^{r_2-2}(f(x)-f(y)) & \left(\frac{ g(x)^{r_1}}{f(x)^{r_1-1}}-\frac{ g(y)^{r_1}}{f(y)^{r_1-1}} \right) \\
        &  \le |g(x)-g(y)|^{r_2-2}(g(x)-g(y)) \left(\frac{g(x)^{r_1-r_2+1}}{f(x)^{r_1-r_2}} - \frac{g(y)^{r_1 - r_2 + 1}}{f(y)^{r_1 - r_2}} \right).
    \end{align*}
    \item[(iii)] Let $\al, \be \ge 1$. Then for $x,y \in \rd$,
    \begin{align*}
        |f(x)-f(y)|^{r_1-2}(f(x)-f(y)) & \left(\frac{ g(x)^{r_1}}{\al f(x)^{r_1-1} + \be f(x)^{r_2-1}}-\frac{g(y)^{r_1}}{\al f(y)^{r_1-1}+\be f(y)^{r_2-1}} \right) \\
        &\le \abs{g(x) -g(y)}^{r_1}.
    \end{align*}
    \item[(iv)] Let $\al, \be \ge 1$. Then for $x,y \in \rd$,
    \begin{align*}
        |f(x)-f(y)|^{r_2-2}(f(x)-f(y)) & \left(\frac{ g(x)^{r_1}}{\al f(x)^{r_1-1} + \be f(x)^{r_2-1}}-\frac{g(y)^{r_1}}{\al f(y)^{r_1-1}+\be f(y)^{r_2-1}} \right) \\
        & \le  \abs{f(x)^{\frac{r_1}{r_2}}-f(y)^{\frac{r_1}{r_2}}}^{r_2}.
    \end{align*}
\end{enumerate}
Moreover, the equality holds in the above inequalities if and only if $f=cg$ a.e. in $\rd$ for some $c \in \R$.
\end{lemma}

\begin{proof}
For the proof of (i), see \cite[Proposition 4.2]{BrFr}. Proof of (ii), (iii), and (iv) follows from \cite[Theorem 2.3 and Remark 2.6]{GiGoMo}.
\end{proof}

The following lemma verifies that certain functions are in the fractional Sobolev space, which we require in the subsequent sections.  

\begin{lemma}\label{test functions are in space}
Let $s \in (0,1)$ and $r_1, r_2 \in (1, \infty)$. Let $u\in W^{s,r_1}_{0}(\Omega)$ be a non-negative function. For $v \in W^{s,r_1}_{0}(\Omega) \cap L^{\infty}(\Om)$, the following functions
$$
\phi_k := \frac{|v|^{r_1}}{u_k^{r_1-1}+u_k^{r_2-1}},\,\psi_k:=\frac{|v|^{r_1}}{u_k^{r_2-1}}, \, \text{ and } \, \eta_k := \frac{|v|^{r_1-r_2+1}}{u_k^{r_1-r_2}} \text{with } r_2<r_1
$$
lie in $W^{s,r_1}_{0}(\Omega)$.
\end{lemma}
\begin{proof}
We only prove that $\phi_k \in W^{s,r_1}_{0}(\Omega)$. For other functions, the proof follows using similar arguments. Clearly, $\phi_k$ is in $L^r_1(\Omega)$ and $\phi_k=0$ in $\Omega^c$, for every $k.$ Next, claim that $[\phi_k]_{s,r_1} < \infty$. In order to show this, for $x,\,y \in \rd$, we calculate
\begin{align*}
& |\phi_k(x)-\phi_k(y)| \\ 
& = \left|\frac{|v(x)|^{r_1}}{u_{k}(x)^{r_1-1}+u_k(x)^{r_2-1}}-\frac{|v(y)|^{r_1}}{u_{k}(y)^{r_1-1}+u_k(y)^{r_2-1}}\right|  \\
& = \left|\frac{|v(x)|^{r_1} - |v(y)|^{r_1} }{u_{k}(x)^{r_1-1}+u_k(x)^{r_2-1}} + \frac{|v(y)|^{ r_1 }\left(u_k(y)^{r_1-1}+u_k(y)^{r_2-1} - (u_k(x)^{r_1-1}+u_k(x)^{r_2-1}\right)}{\left(u_{k}(x)^{r_1-1}+u_k(x)^{r_2-1}\right)\left(u_k(y)^{r_1-1} + u_k(y)^{r_2-1}\right)}\right|  \\
& \leq \left(k^{r_1-1} + k^{r_2-1}\right) \left||v(x)|^{r_1} - |v(y)|^{ r_1 }\right| \\
& \quad + \norm{v}^{r_1}_{\infty} \frac{\left|u_k(y)^{r_1-1}-u_k(x)^{r_1-1}\right| + \left|u_k(y)^{r_2-1}-u_k(x)^{r_2-1}\right|}{\left(u_{k}(x)^{r_1-1}+u_k(x)^{r_2-1}\right)\left(u_k(y)^{r_1-1}+u_k(y)^{r_2-1}\right)}.
\end{align*}
Using (iii) of Lemma \ref{inequality}, we get 
\begin{equation*}
    \begin{split}
        |\phi_k(x)-\phi_k(y)| & \le r_1\left(k^{r_1-1}+k^{r_2-1}\right)\left(|v(x)|^{r_1-1}+|v(y)|^{r_1-1}\right)|v(x)-v(y)| \\
        & \quad + (r_1-1)\norm{v}^{r_1}_{\infty}\frac{\left(u_{k}(x)^{r_1-2}+u_k(y)^{r_1-2}\right)}{u_{k}(x)^{r_1-1}u_k(y)^{r_1-1}}|u_k(x)-u_k(y)|
         \\
        & \quad + (r_2-1)\norm{v}^{r_1}_{\infty}\frac{\left(u_{k}(x)^{r_2-2}+u_k(y)^{r_2-2}\right)}{u_{k}(x)^{r_2-1}u_k(y)^{r_2-1}}|u_k(x)-u_k(y)|.
    \end{split}
\end{equation*}
Now using $u_k^{-1} \le k$ and $v \in L^{\infty}(\Om)$, there exists $C=C(r_1,r_2,k,\norm{v}_{\infty})$ such that
\begin{align*}
    |\phi_k(x)-\phi_k(y)| & \le C\left(|v(x)-v(y)|+|u_k(x)-u_k(y)|\right) = C\left(|v(x)-v(y)|+|u(x)-u(y)|\right).
\end{align*}
Therefore, $\phi_k \in W^{s,r_1}_{0}(\Omega)$ follows as $[v]_{s,r_1},\,[u]_{s,r_1}<\infty.$ 
This completes the proof.
\end{proof}

\section{Linear independence of the first eigenfunctions}\label{section5}
This section is devoted to proving the linear independency of the first Dirichlet eigenfunctions of the fractional $p$-Laplacian and the fractional $q$-Laplacian. Throughout the section, we assume that $\Omega \subset \rd$ is a bounded open set of class $\mathcal{C}^{1,1}$. For brevity, we denote the first eigenpair of \eqref{nonlocal Dirichlet EV} by $(\la^1_{s,r},u)$. From Proposition \ref{first eigenvalue}, $u>0$ in $\Omega$, $u=0$ in $\rd\setminus\Omega$ and $u \in \C(\overline{\Om})$. Therefore, $u$ attains its maximum in ${\Om}$. Due to the translation invariance of the fractional $r$-Laplacian, we assume that $\Om$ contains the origin and the maximum point for $u$ is the origin. Now for $\tau>0$, we consider $\Om_\tau := \{z \in \rd : \tau z \in \Om \}$ and define $u_{\tau} : \rd \ra \R$ as follows:
\begin{equation*}
   u_\tau(x) := \left\{\begin{array}{ll} 
            \frac{ u(0)-u(\tau x)}{\tau^{sr^{\prime}}} , & \text {for }  x\in\Om_\tau; \\ 
            \frac{ u(0)}{\tau^{sr^{\prime}}}, & \text{for } \; x \in \rd \setminus \Om_\tau. \\
             \end{array} \right.
\end{equation*}
The following result demonstrates a property of the above function, which plays an essential role in proving \eqref{LI}.
\begin{lemma}[Blow-up lemma]\label{blowup lemma}
Let $r \in (1,  \infty)$ and $s \in (0,1)$. If $\tau_n \to 0$, as $n \to \infty$, then there exists a subsequence denoted by $(\tau_n)$ such that $u_{\tau_n} \to \tilde{u}$ in $\C_{\text{loc}}(\rd)$ as $n \to \infty$. Moreover, $\tilde{u} \in W^{s,r}_{\text{loc}}(\rd)\cap \C(\rd)$ is non-negative, and satisfies the following equation weakly:
\begin{equation}\label{u Bar eq}
(-\Delta)^s_r v = -\lambda^1_{s,r} u(0)^{r-1} \text{ in }\rd,
\end{equation}
and $\tilde{u}(0) = 0.$
\end{lemma}

\begin{proof}
Note that for any $\tau > 0$, $u_\tau \geq 0$, since $u(0)$ is the maximum value for $u$ in $\overline{\Omega}$. Using the fact that $(\la^1_{s,r}, u(\tau x))$ is the first eigenpair for fractional $r$-Laplacian on $\Omega_\tau$, we obtain that the following equation holds weakly:
\begin{equation}\label{u_tau eq}
  (-\De)^s_r u_\tau(x) = -(-\De)^s_r u(\tau x) = -\lambda^1_{s,r}u(\tau x)^{r-1} \text{ in }\Om_\tau, \quad u_\tau = \frac{ u(0)}{\tau^{sr^{\prime}}} \text{ in } \rd \setminus \Om_\tau. 
 \end{equation}
For each $\tau>0$, using Proposition \ref{first eigenvalue} and \cite[Theorem~3.13]{BrPa}, we get $u_{\tau} \in \C(\Om_{\tau})$.  Now we divide our proof into two steps. In the first step, we show that $u_{\tau_n}\to\tilde{u}$ in $\C_{\text{loc}}(\rd)$ as $n\to\infty$. In the second step, we prove $\tilde{u}$ is a weak solution to 
\eqref{u Bar eq}.
\smallskip

\noi \textbf{Step 1:} Take a ball $B_R(0)$ such that $\overline{B_{4R}(0)}\subset\Omega_\tau$. We choose $\sigma_1>0$ as follows
\begin{align*}
    \sigma_1:=\frac{d\ga}{s}, \text{ where } \ga >1.
\end{align*}
By the nonlocal Harnack inequality (see \cite[Theorem 2.2]{GiKuSr}), there exists $C=C(d,s,r)$ such that
\begin{align}\label{uniform bound 0}
    \max_{\substack{B_R(0)}}u_\tau \leq C\left(\min_{\substack{B_{2R}(0)}}u_\tau+\norm{\lambda^1_{s,r}\,u^{r-1}}_{L^{\sigma_1}(B_{2R}(0))}^{\frac{1}{r-1}}\right) = C \norm{\lambda^1_{s,r}\,u^{r-1}}_{L^{\sigma_1}(B_{2R}(0))}^{\frac{1}{r-1}},
\end{align}
In \eqref{uniform bound 0} the last equality follows from the fact $\min_{\substack{B_{2R}(0)}}u_\tau=0$, because origin is the maximum point of $u$ in $\Om$. Further, for $r\geq 2$ we immediately get $(r-1) \sigma_1 > 1$, and for $1<r<2$ we choose 
\begin{equation*}
    \ga> \left\{\begin{array}{ll} 
            \text{max}\left\{\frac{sr}{dr+2(sr-d)},1\right\} , & \text {if }  sr<d; \\ 
            \text{max}\left\{\frac{s}{d(r-1)},1\right\}, & \text{if} \; sr \ge d,  \\
             \end{array} \right.
\end{equation*}
to get $(r-1)\sigma_1>1$. Then Proposition \ref{first eigenvalue}-(v) and \eqref{uniform bound 0} yield
\begin{align}\label{uniform bound 1}
    \max_{\substack{B_R(0)}}u_\tau \leq C \norm{\lambda^1_{s,r}\,u^{r-1}}_{L^{\sigma_1}( \rd )}^{\frac{1}{r-1}} \le C,
\end{align}
where $C=C(d, s, r, \lambda^1_{s,r},\norm{u}_{L^{(r-1)\sigma_1}(\rd)})$.
Next, we define the following exponent
\begin{align}\label{range}
  \Theta(d,s,r,\sigma_1) := \min \left\{ \frac{1}{r-1} \left( sr - \frac{d}{\sigma_1} \right), 1 \right\}.
\end{align}
Then, applying the regularity estimate \cite[Theorem~1.4]{BrLiSc} when $r \ge 2$ and \cite[Theorem~1.2]{GaLi} when $1<r < 2$,  for the problem \eqref{u_tau eq} we get the following H\"{o}lder regularity estimate of the weak solution $u_{\tau}$ for any $s<\delta<\Theta(d,s,r,\sigma_1)$:
\begin{equation}\label{Holder norm est}
    \begin{split}
        [u_\tau]_{C^{0,\delta}(B_{R/8}(0))}
        &\leq\frac{C}{R^{\delta}} \left[\max_{\substack{B_R(0)}}u_\tau + \Tail_r(u_\tau;0,R) + \left(R^{sr - \frac{d}{\sigma_1} } \la^1_{s,r}\norm{u^{r-1}}_{L^{\sigma_1} (B_R(0))} \right)^{\frac{1}{r-1}} \right]
         \\
         & :=\frac{C}{R^\delta}\left[\max_{\substack{B_R(0)}}u_\tau+I_1+I_2^{\frac{1}{r-1}}\right],
    \end{split}
\end{equation}
where $C=C(d,s,r)$. We now estimate the last two terms $I_1,\,I_2$ of \eqref{Holder norm est} as follows:
\\
\smallskip
\noindent\textbf{Estimate of $I_2$:} Choose $a>0$ such that $a >\ga r-1$. Then, by the change of variable we have
\begin{equation}\label{I_2 est blowup lemma}
    \begin{split}
        I_2:=\lambda^1_{s,r} R^{ s(r - \frac{1}{\ga}) } \left( \int_{B_R(0)} \abs{u(y)}^{(r - 1) \sigma_1} \, \dy \right)^{\frac{1}{\sigma_1}} & =  \lambda^1_{s,r} \frac{R^{ s(r - \frac{1}{\ga}) }}{R^{ \frac{s a}{\ga} }} \left( \int_{B_{R^{ a + 1}}(0)} \abs{u(z)}^{ (r-1) \sigma_1} \, \dz \right)^{\frac{1}{\sigma_1}} \\
        & \le \lambda^1_{s,r} R^{ s(r - \frac{1}{\ga}) - \frac{s a}{\ga}}  \left( \int_{\rd} \abs{u(z)}^{ (r-1) \sigma_1} \, \dz \right)^{\frac{1}{\sigma_1}},
    \end{split}
\end{equation}
 where we see that $r - \frac{1}{\ga} < \frac{a}{\ga}$. 
\\
\smallskip
\noindent\textbf{Estimate of $I_1$:}
Note that
\begin{equation}\label{relation of tails}
I_1:=\Tail_r(u_\tau;0,R)\leq C(r)\left(\Tail_r(u_{\tau}^+;0,R)+\Tail_r(u_{\tau}^-;0,R)\right)=C(r)\Tail_r(u_{\tau}^+;0,R),
\end{equation}
where the last equality follows from the non-negativity of $u_\tau$. To estimate $\Tail_r(u_{\tau}^+;0,R)$, let $R_1 = 4R$ and $\ell:=\max_{\substack{B_{R}(0)}}u_\tau$. Take $ \phi \in \cc(B_{R})$ satisfying $0\leq\phi\leq1$, $\phi=1$ in $B_{\frac{R}{2}}$ and $\abs{\nabla\phi}\leq \frac{8}{R}$. We use the test function $\eta :=(u_\tau-2\ell)\phi^p$ in the weak formulation of $u_\tau$ and then proceed similarly as in \cite[Lemma~4.2]{CaKuPa}) to get a constant $C=C(d,s,r)$ such that 
\begin{equation*}
\begin{split}
C\ell\abs{B_{R}}R^{-sr}\Tail_r(u_{\tau}^+;0,R)^{r-1}
& \leq C\ell^{r}\abs{B_{R}}R^{-sr}+\lambda^1_{s,r}\int_{B_{R}}u^{r-1}\eta \, \dx \\
&\leq C\ell^{r}\abs{B_{R}}R^{-sr}+3\ell\lambda^1_{s,r}\int_{B_{R}}u^{r-1} \, \dx
\\
&\leq C\ell^{r}\abs{B_{R}}R^{-sr}+3\ell\lambda^1_{s,r}\abs{B_{R}}^{\frac{1}{ \sigma_1'}} \norm{u^{r-1}}_{L^{\sigma_1}(B_R)},
\end{split}
\end{equation*}
where in the above estimates we used the fact $|u-2\ell|\leq3\ell$ in $B_{R}$.
This implies that
\begin{equation}\label{tail plus upper bdd}
\begin{split}
    \Tail_r(u_{\tau}^+;0,R)
   & \leq C\left(\max_{\substack{B_{R}(0)}}u_\tau+\left(R^{sr - \frac{d}{\sigma_1} }\lambda^1_{s,r}\norm{u^{r-1}}_{L^{\sigma_1}(B_R)}\right)^{\frac{1}{r-1}}\right) \\
   &\le C\left(\max_{\substack{B_{R}(0)}}u_\tau + I_2^{\frac{1}{r-1}}\right).
\end{split}
\end{equation}
Now, plugging the estimates \eqref{uniform bound 1}, \eqref{I_2 est blowup lemma}, \eqref{relation of tails}, \eqref{tail plus upper bdd} into \eqref{Holder norm est}, we thus obtain 
\begin{equation}\label{final Holder norm est-2}
[u_\tau]_{\C^{0, \delta}(B_{R/8}(0))}\leq\frac{C}{R^{\delta+\epsilon}},
\end{equation}
where $C=C(d,s,r,\lambda^1_{s,r}, \norm{u}_{L^{ (r-1) \sigma_1}(\rd)})$ and $\epsilon:=\frac{1}{r-1}\left(\frac{1+a-r\ga}{\ga}\right) > 0$. Let $K \subset \rd$ be any compact set. Observe that $\Omega_{\tau}$ becoming $\rd$ when $\tau$ is sufficiently small. Thus, we can choose $R>1$ and $0<\tau_0<<1$ such that $K \subset B_{\frac{R}{8}}(0) \subset \Om_{\tau}$ for all $\tau \in (0, \tau_0)$. Therefore, we use \eqref{uniform bound 1} and \eqref{final Holder norm est-2} to obtain the following uniform estimate for all $\tau \in (0,\tau_0)$:
\begin{equation}\label{C_alpha bdd}
  \max_{K} u_{\tau} \le C, \text{ and }  [u_\tau]_{\C^{0, \delta }(K)} \le C,
\end{equation}
where $C$ is independent of both $\tau$ and $K$. Next, for a sequence $(\tau_n)$ converging to zero,  we consider the corresponding sequence of functions $(u_{\tau_n})$. Using \eqref{C_alpha bdd} we can show that $(u_{\tau_n})$ is equicontinuous and uniformly bounded in $K$. Therefore, applying the Arzela-Ascoli theorem, up to a subsequence, $u_{\tau_n}\to\tilde{u}\text{ in }\C(K)$. Thus we have 
\begin{equation}\label{C_loc convergence}
   u_{\tau_n}\to\tilde{u}\text{ in } \C_{\text{loc}}(\rd), \text{ as } n \to \infty. 
\end{equation}

\noi \textbf{Step 2:} Recalling the weak formulation of \eqref{u_tau eq} for $\tau>0$ be any, 
\begin{equation}\label{weak formulation of u_tau}
    \begin{split}
        \iint\limits_{\rd\times\rd}\frac{|u_\tau(x)-u_\tau(y)|^{r-2}(u_\tau(x)-u_\tau(y))(\phi(x)-\phi(y))}{|x-y|^{d+sr}} \,  \dx \dy \\ = -\la^1_{s,r} \int_{\Omega_\tau}  u(\tau x)^{r-1} \phi(x) \, \dx, \quad \forall \, \phi\in\cc(\Omega_\tau).
    \end{split}
\end{equation}
Let $v\in\cc(\rd)$ and let $\text{supp}(v) := K$. Since $\Omega_{\tau_n}$ is becoming $\rd$, as $\tau_n\to0$, there exists $n_0 \in \mathbb{N}$ such that $K \subset\Om_{\tau_n}$ for all $n\geq n_0.$ Hence, from \eqref{weak formulation of u_tau} for every $n \ge n_0$, we write
\begin{equation}\label{weak formulation of u_taun}
    \begin{split}
\iint\limits_{\rd\times\rd}\frac{|u_{\tau_n}(x)-u_{\tau_n}(y)|^{r-2}(u_{\tau_n}(x)-u_{\tau_n}(y))(v(x)-v(y))}{|x-y|^{d+sr}} \, & \dx \dy  
     \\
     = -\la^1_{s,r} \int_{K}  u(\tau_n x)^{r-1} v(x) \, \dx.
    \end{split}
\end{equation}
We pass the limit as $n\to\infty$ in the R.H.S of \eqref{weak formulation of u_taun}, to get
\begin{align}\label{right hand side}
\lim\limits_{n\to\infty} \int_{K}  u(\tau_n x)^{r-1}v(x) \, \dx = \lim\limits_{n\to\infty}\int_{\rd}u(\tau_n x)^{r-1}v(x)\chi_{K}(x) \, \dx  = \int_{K} u(0)^{r-1}v(x)\,\dx,
\end{align}
where the last equality in \eqref{right hand side} follows using the dominated convergence theorem. Again, applying the dominated convergence theorem, we have
\begin{align*}
    \text{L.H.S of \eqref{weak formulation of u_taun}} & = \lim\limits_{k\to\infty}\II{B_k(0)\times B_k(0)}\frac{|u_{\tau_n}(x)-u_{\tau_n}(y)|^{r-2}(u_{\tau_n}(x)-u_{\tau_n}(y))(v(x)-v(y))}{|x-y|^{d+sr}} \,  \dx \dy
    \\
    & := \lim\limits_{k\to\infty}\II{B_k(0)\times B_k(0)} F_n(x,y)\,\dx\dy.
\end{align*}
\noi\textbf{Claim:} Now, we establish 
\begin{align}\label{limit intercahnge}
\lim\limits_{n\to\infty}\lim\limits_{k\to\infty}\II{B_k(0)\times B_k(0)} F_n(x,y)\,\dx\dy
&=\lim\limits_{k\to\infty}\lim\limits_{n\to\infty}\II{B_k(0)\times B_k(0)} F_n(x,y)\,\dx\dy.
\end{align}
\noi Proof of Claim: To show \eqref{limit intercahnge}, for any fixed $k \in \mathbb{N}$ we first prove that 
$$
F_n(x,y)\xrightarrow{n\to\infty}F(x,y):=\frac{|\tilde{u}(x)-\tilde{u}(y)|^{r-2}(\tilde{u}(x)-\tilde{u}(y))(v(x)-v(y))}{|x-y|^{d+sr}} \, \text{ in }L^1\left(B_k(0)\times B_k(0)\right).
$$
It is easy to see from \eqref{C_loc convergence} that $F_n(x,y)\xrightarrow{n\to\infty}F(x,y)$ pointwise. Now for $x,y\in B_k(0)$, and using the uniform boundedness of $(u_{\tau_n})$ (see \eqref{C_alpha bdd}), we have
\begin{align*}
    |F_n(x,y)| =\frac{|u_{\tau_n}(x)-u_{\tau_n}(y)|^{r-1}|v(x)-v(y)|}{|x-y|^{d+sr}}
    & \leq [u_{\tau_n}]^{r-1}_{\C^{0, \delta}\left(\overline{B_k(0)}\right)}\frac{|v(x)-v(y)|}{|x-y|^{d + sr -\delta(r-1)}}
    \\
    & \le C\,\frac{|v(x)-v(y)|}{|x-y|^{d + sr-\delta(r-1)}}, 
\end{align*}
where the constant $C$ does not depend on $n$. By Fubini's theorem, we get for any fixed $k\in\mathbb{N}$
\begin{equation*}
    \begin{split}
        \II{B_k(0)\times B_k(0)}\frac{|v(x)-v(y)|}{|x-y|^{d+sr-\delta(r-1)}}\,\dx \dy
      & \leq\II{B_k(0)\times B_{2k}(0)}\frac{|v(x)-v(x+z)|}{|z|^{d+sr-\delta(r-1)}}\,\dz \dx
      \\
    & =\II{B_k(0)\times B_{2k}(0)}\left(\int_{0}^1\frac{|\nabla v(x+tz)|}{|z|^{d+sr-\delta(r-1)-1}} \dt\right)\dz\dx \\
    & \leq\int_{B_{2k}(0)}\int_0^1\frac{\norm{\nabla v}_{L^1(\rd)}}{|z|^{d+sr-\delta(r-1)-1}}\dt\dz
    \\
    & =C\norm{\nabla v}_{L^1(\rd)}<\infty,\, \text{ since } \delta > \frac{sr-1}{r-1}.
    \end{split}
\end{equation*}
Thus, applying the dominated convergence theorem, we conclude $F_n \xrightarrow{n\to\infty} F$ in $L^1\left(B_k(0)\times B_k(0)\right)$.
Also, it is easy to verify that for any fixed $n\in\mathbb{N}$
\begin{equation*}
    \II{\rd\times\rd}F_{n}(x,y)\chi_{B_k(0)}(x)\chi_{B_k(0)}(y)\, \dx\dy\xrightarrow{k\to\infty} \II{\rd\times\rd}F_n(x,y)\,\dx\dy.
\end{equation*}
Again, from the Fatou's lemma, \eqref{weak formulation of u_taun}, and \eqref{right hand side} we get 
\begin{align*}
   \II{\rd\times\rd}F(x,y)\,\dx\dy \le \lowlim_{n \ra \infty}  \II{\rd\times\rd}F_n(x,y)\,\dx\dy =- \la^1_{s,r} \lowlim_{n \ra \infty} \int_{K}  u(\tau_n x)^{r-1}v(x) \, \dx \le C.
\end{align*}
Next, for $n, k \in \mathbb{N}$, we consider the double sequence of functions $(F_{n,k})$ defined as
\begin{equation*}
        F_{n,k}(x,y) := F_{n}(x,y)\chi_{B_k(0)}(x)\chi_{B_k(0)}(y), \text{ for }  x, y \in \rd.
\end{equation*}
We claim that
\begin{equation}\label{double limit}
\lim\limits_{\substack{n\to\infty\\k\to\infty}}\,\,\II{\rd\times\rd}F_{n,k}(x,y)\,\dx\dy=\II{\rd\times\rd}F(x,y)\,\dx\dy.
\end{equation}
Again, using \eqref{C_loc convergence}, $F_{n,k}(x,y)\xrightarrow{n,\,k\to\infty}F(x,y)$ pointwise a.e. in $\rd$. Further, for $x,\,y\in\rd$, using the uniform estimate \eqref{C_alpha bdd} we have
\begin{align*}
    |F_{n,k}(x,y)| & =|F_n(x,y)|\chi_{B_k(0)}(x)\chi_{B_k(0)}(y) \\
    & =\frac{|u_{\tau_n}(x)-u_{\tau_n}(y)|^{r-1}|v(x)-v(y)|}{|x-y|^{d+sr}} \chi_{B_k(0)}(x)\chi_{B_k(0)}(y)\\
    & \leq [u_{\tau_n}]^{r-1}_{\C^{0, \delta} \left(\overline{B_k(0)}\right)}\frac{|v(x)-v(y)|}{|x-y|^{d + sr-\delta(r-1)}}\chi_{B_k(0)}(x)\chi_{B_k(0)}(y)
    \\
    & \le C \frac{|v(x)-v(y)|}{|x-y|^{ d + sr-\delta(r-1)}}\chi_{B_k(0)}(x)\chi_{B_k(0)}(y), 
\end{align*}
where the constant $C$ is independent of both $n$ and $k.$ Moreover, from  the fact that $\delta > \frac{sr-1}{r-1}$, $$\II{\rd\times\rd}\frac{|v(x)-v(y)|}{|x-y|^{d + sr-\delta(r-1)}}\,\dx \dy < \infty,$$ if we choose $\delta < \frac{sr}{r-1}$.
Thus, \eqref{double limit} follows by again using the dominated convergence theorem. Hence, by the standard result for interchanging double limits, we obtain \eqref{limit intercahnge}. 

\noi Therefore, taking the limit as $n\to\infty$ in the L.H.S of  \eqref{weak formulation of u_taun} and using \eqref{limit intercahnge} we obtain
\begin{equation}\label{double limits}
    \begin{split}
        \lim\limits_{n\to\infty}\lim\limits_{k\to\infty}\II{B_k(0)\times B_k(0)} F_n(x,y)\,\dx\dy
        &=\lim\limits_{k\to\infty}\lim\limits_{n\to\infty}\II{B_k(0)\times B_k(0)} F_n(x,y)\,\dx\dy
        \\
        &=\lim\limits_{k\to\infty}\II{B_k(0)\times B_k(0)} F(x,y)\,\dx\dy \\
        &=\II{\rd\times\rd}F(x,y)\,\dx\dy,
    \end{split}
\end{equation}
Thus, using \eqref{weak formulation of u_taun}, \eqref{right hand side}, and \eqref{double limits} we get
\begin{equation*}
    \begin{split}
        \II{\rd\times\rd} \frac{|\tilde{u}(x)-\tilde{u}(y)|^{r-2}(\tilde{u}(x)-\tilde{u}(y))(v(x)-v(y))}{|x-y|^{d+sr}} \, \dx\dy \\
    =  -\la^1_{s,r}\int_{\rd}u(0)^{r-1}v(x)\,\dx, \quad \forall \, v \in \cc(\rd).
    \end{split}
\end{equation*}
 Moreover, we also have $\tilde{u}\in W^{s,r}_{\text{loc}}(\rd)\cap \C(\rd)$ provided $s<\delta$. Hence, $\tilde{u}$ is a weak solution of \eqref{u Bar eq}. Again, since $u_{\tau_n}\geq0$, $u_{\tau_n}(0)=0$, from \eqref{C_loc convergence} we arrive at $\tilde{u}\geq0$ with $\tilde{u}(0)=0.$  This completes the proof of the lemma. 
\end{proof}

\noi \textbf{Proof of Theorem \ref{linearindependent}}: For simplicity of notation, we denote $u_0 = \phi_{s_1,p}$ and $v_0 = \phi_{s_2, q}$.
We argue by contradiction. Suppose $u_0 = cv_0$ for some non-zero $c \in \R$. Without loss of any generality, we can assume that $u_0=v_0.$ By Proposition \ref{first eigenvalue}, $u_0$ is bounded, $u_0>0$ in $\Omega$ and is in $\C( \overline{ \Omega })$. This guarantees that $u_0$ has a global extremum point.  Since the operator $(-\Delta)_p^{s_1}$ is translation invariant, we can assume that the origin is such a point. Now for $\tau>0$, define
\begin{equation}\label{u_0 tau eq}
   u_\tau(x) := \left\{\begin{array}{ll} 
            \frac{ u_0(0)-u_0(\tau x)}{\tau^{s_1\p}} , & \text {for }  x\in\Om_\tau; \vspace{0.1 cm} \\ 
            \frac{ u_0(0)}{\tau^{s_1\p}}, & \text{for } \; x \in \rd \setminus \Om_\tau, \\
             \end{array} \right.
\end{equation}
where $\Omega_\tau:=\{x\in\rd:\tau x\in\Omega\}$.
Then by Blow-up lemma \ref{blowup lemma}, there exists a sequence $\tau_n\to 0$ such that $u_{\tau_n}\to\tilde{u}$ in $\C_{\text{loc}}(\rd)$, where $\tilde{u}$ is a non-negative solution of 
\begin{equation}\label{u_0 bar eq}
(-\Delta)^{s_1}_p v=-\lambda^1_{s_1,p}u_0(0)^{p-1} \text{ in }\rd,
\end{equation}
and $\tilde{u}(0)=0.$ Again, by the change of variable we deduce
\begin{align*}
    (-\Delta)_q^{s_2}u_\tau(x) & =\text{P.V.}\int_{\rd}\frac{|u_\tau(x)-u_\tau(y)|^{q-2}(u_\tau(x)-u_\tau(y))}{|x-y|^{d+s_2q}} \,\dy
    \\
    & =-\frac{1}{\tau^{s_1\p(q-1)}}\text{P.V.}\int_{\rd}\frac{|u_0(\tau x)-u_0(\tau y)|^{q-2}(u_0(\tau x)-u_0(\tau y))}{|x-y|^{d+s_2q}} \,\dy
    \\
    & =-\tau^{s_2q-s_1\p(q-1)}\,\text{P.V.}\int_{\rd}\frac{|u_0(\tau x)-u_0(y)|^{q-2}(u_0(\tau x)-u_0(y))}{|\tau x-y|^{d+s_2q}} \,\dy
    \\
    & =-\tau^{s_2q-s_1\p(q-1)}(-\Delta)_q^{s_2}u_0(\tau x)=-\tau^{s_2q-s_1\p(q-1)}\la^1_{s_2,q}\,u_0(\tau x)^{q-1}.
\end{align*}
This implies that for each $\tau>0$, $u_\tau$ given by \eqref{u_0 tau eq} satisfies the following equation weakly
$$
(-\Delta)_q^{s_2}v=-\tau^{s_2q-s_1\p(q-1)}\la^1_{s_2,q}\,u_0(\tau x)^{q-1} \text{ in }\Omega_\tau.
$$
Using $\frac{s_1\p}{\q}<s_2$ we again proceed as in Blow-up lemma \ref{blowup lemma}, to obtain that $\tilde{u} \ge 0$ is also a weak solution of the following equation:
\begin{equation*}
    (-\Delta)_q^{s_2} v=0\,\text{ in }\rd.
\end{equation*}
Therefore, by the strong maximum principle \cite[Theorem~1.4]{PeQu}, we conclude $\tilde{u}=0$ in $\rd,$ which gives a contradiction to \eqref{u_0 bar eq} as $u_0(0)>0$. Thus, the set $\{u_0,v_0\}$ is linearly independent.
\qed

\section{$L^\infty$ bound and maximum principle}\label{section6}
In this section, under the presence of multiple exponents $(s_1,p),(s_2,q)$ and parameters $(\al,\be)$, 
we first prove that every nonnegative weak solution of \eqref{main problem} is bounded in $\rd$. Afterwards, we state a strong maximum principle. 

\begin{theorem}[Global $L^\infty$ bound]\label{bounded solution}
Let $0<s_2<s_1<1<q<p<\infty$ and let $\Omega\subset\rd$ be a bounded open set. Assume that $u\in\wpso$ is a nonnegative solution of \eqref{main problem}. Then $u\in L^{\infty}(\rd)$.
\end{theorem}
\begin{proof}
\noi $\bm{d>s_1p\,}$\textbf{:} Let $M\geq0$, define $u_M=\min\{u,M\}$. Clearly $u_M$ is non-negative and is in $L^{\infty}(\Om)$. Since $u\in\wpso$, then $u_M\in\wpso.$ Fixed $\sigma\geq1$, define $\phi=u_M^{\sigma}$. Then, $\phi\in\wpso.$ Thus taking $\phi$ as a test function in the weak formulation of $u$, we have
\begin{equation}\label{Appendix: solution bdd eq1}
    \begin{split}
        & \iint\limits_{\rd\times\rd}|u(x)-u(y)|^{p-2}(u(x)-u(y))(\phi(x)-\phi(y)) \, d\mu_1
       \\
       &+\iint\limits_{\rd\times\rd}|u(x)-u(y)|^{q-2}(u(x)-u(y))(\phi(x)-\phi(y))\, d\mu_2
       \\
       &=\al\intom u(x)^{p-1}\phi(x) \,\dx+\be\intom u(x)^{q-1}\phi(x)\,\dx \leq \al\intom u(x)^{p+\sigma-1}\,\dx+\be\intom u(x)^{q+\sigma-1}\,\dx.
    \end{split}
\end{equation}
Now, using \cite[Lemma~C.2]{BrLiPa} we estimate
\begin{equation*}
    \begin{split}
        I_1 &  :=\iint\limits_{\rd\times\rd}|u(x)-u(y)|^{p-2}(u(x)-u(y))(\phi(x)-\phi(y)) \, d\mu_1\\
        & \geq\frac{\sigma p^p}{(\sigma+p-1)^p}\iint\limits_{\rd\times\rd}\left|u_M(x)^{\frac{\sigma+p-1}{p}}-u_M(y)^{\frac{\sigma+p-1}{p}}\right|^{p}\,d\mu_1
       \\ & \geq \frac{C(d,s_1,p)\sigma p^p}{(\sigma+p-1)^p}
       \left(\int_{\rd}\left(u_M(x)^{\frac{\sigma+p-1}{p}} \right)^{p^*_{s_1}}\,\dx\right)^{\frac{p}{p^*_{s_1}}},
    \end{split}
\end{equation*}
where in the last inequality we use $\wpso \hookrightarrow L^{p^*_{s_1}}(\rd)$. Since $s_2q < d$, using $\wqst \hookrightarrow L^{q^*_{s_2}}(\rd)$ we estimate $I_2$ as 
\begin{align*}
I_2 & :=\iint\limits_{\rd\times\rd}|u(x)-u(y)|^{q-2}(u(x)-u(y))(\phi(x)-\phi(y))\, d\mu_2\\
&\geq \frac{C(d,s_2,q)\sigma q^q}{(\sigma+q-1)^q}
\left(\int_{\rd}\left(u_M(x)^{\frac{\sigma+q-1}{q}}\right)^{q^*_{s_2}}\,\dx\right)^{\frac{q}{q^*_{s_2}}}.
\end{align*}
Plugging the estimates of $I_1$ and $I_2$ into \eqref{Appendix: solution bdd eq1} we obtain
\begin{equation*}
    \begin{split}
        \frac{C(d,s_1,p)\sigma p^p}{(\sigma+p-1)^p}
\left(\int_{\rd}\left(u_M(x)^{\frac{\sigma+p-1}{p}} \right)^{p^*_{s_1}}\,\dx\right)^{\frac{p}{p^*_{s_1}}}
+\frac{C(d,s_2,q)\sigma q^q}{(\sigma+q-1)^q}
\left(\int_{\rd}\left(u_M(x)^{\frac{\sigma+q-1}{q}} \right)^{q^*_{s_2}}\,\dx\right)^{\frac{q}{q^*_{s_2}}}
\\
\leq\al\intom u(x)^{p+\sigma-1}\,\dx+\be\intom u(x)^{q+\sigma-1}\,\dx.
    \end{split}
\end{equation*}
Letting $M\to\infty$ in above, the monotone convergence theorem yields
\begin{equation}\label{M to infty}
    \begin{split}
        \frac{C(d,s_1,p)\sigma p^p}{(\sigma+p-1)^p}
\left(\int_{\rd}\left(u(x)^{\frac{\sigma+p-1}{p}}\right)^{p^*_{s_1}}\,\dx\right)^{\frac{p}{p^*_{s_1}}}
+\frac{C(d,s_2,q)\sigma q^q}{(\sigma+q-1)^q}
\left(\int_{\rd}\left(u(x)^{\frac{\sigma+q-1}{q}}\right)^{q^*_{s_2}}\,\dx\right)^{\frac{q}{q^*_{s_2}}}
\\
\leq\al\intom u(x)^{p+\sigma-1}\,\dx+\be\intom u(x)^{q+\sigma-1}\,\dx.
    \end{split}
\end{equation}
\noi \textbf{Claim:} For $\sigma_1 := p^*_{s_1} - p +1$, $u^{\sigma_1 + p -1} \in L^{\frac{p^*_{s_1}}{p}}(\rd)$. \\
\noi By taking $\sigma = \sigma_1$, we obtain from \eqref{M to infty} that 
\begin{align}\label{claim 1.0}
    \frac{C(d,s_1,p)\sigma p^p}{(p^*_{s_1})^p} \left(\int_{\rd}u(x)^{\frac{p^*_{s_1}}{p}p^*_{s_1}} \,\dx\right)^{\frac{p}{p^*_{s_1}}} \le \al\intom u(x)^{p^*_{s_1}}\,\dx+\be\intom u(x)^{q+\sigma_1 -1}\,\dx.
\end{align}
Notice that $q+ \sigma_1 -1 = q + p^*_{s_1} -p < p^*_{s_1}$ (as $p>q$). Set $a_1:=\frac{p^*_{s_1}}{q+ \sigma_1 -1}$. By applying the H\"{o}lder's inequality with conjugate pair $(a_1, a_1')$ we estimate the second integral of  \eqref{claim 1.0} as
\begin{align}\label{claim 1.1}
    \intom u(x)^{q+\sigma_1 -1}\,\dx \le \left( \intom u(x)^{p^*_{s_1}} \, \dx \right)^{\frac{1}{a_1}} |\Om|^{\frac{1}{a_1'}}.
\end{align}
For $R>1$, consider the set $A:= \{ x \in \Om : u(x) \le R \}$ and $A^c=\Om\setminus A$. We estimate the first integral of the R.H.S of \eqref{claim 1.0} as follows:
\begin{equation}\label{claim 1.2}
    \intom u(x)^{p^*_{s_1}}\,\dx= \left(\int_{A}+\int_{A^c}\right) u(x)^{p^*_{s_1}}\,\dx\le R^{p^*_{s_1}}\abs{\Om} + \abs{A^c}^{\frac{p^*_{s_1}-p}{p^*_{s_1}}} \left( \int_{A^c} u(x)^{p^*_{s_1} \frac{p^*_{s_1}}{p}}\,\dx \right)^{\frac{p}{p^*_{s_1}}}.
\end{equation}
We choose $R>1$ so that 
\begin{align*}
 \alpha \frac{(p^*_{s_1})^p}{C(d,s_1,p)\sigma_1 p^p} \abs{A^c}^{\frac{p^*_{s_1}-p}{p^*_{s_1}}} \le \frac{1}{2}.
\end{align*}
Therefore, combining \eqref{claim 1.0}, \eqref{claim 1.1}, and \eqref{claim 1.2} we obtain 
\begin{align*}
    \frac{1}{2} \left(\int_{\rd}u(x)^{\frac{p^*_{s_1}}{p}p^*_{s_1}} \,\dx\right)^{\frac{p}{p^*_{s_1}}} \le \frac{(p^*_{s_1})^p}{C(d,s_1,p)\sigma_1 p^p} 
 \left( \al \abs{\Om} R^{p^*_{s_1}} + \be |\Om|^{\frac{1}{a_1'}} \left( \intom u(x)^{p^*_{s_1}} \, \dx \right)^{\frac{1}{a_1}}  \right).
\end{align*}
Thus, $u^{\sigma_1 + p -1} \in L^{\frac{p^*_{s_1}}{p}}(\rd)$ for $\sigma_1 := p^*_{s_1} - p +1$. Set $a_2:= \frac{p^*_{s_1}+\sigma - 1}{p+\sigma - 1}$ and $a_3 := \frac{p^*_{s_1}+\sigma - 1}{q+\sigma - 1}$. Using the Young's inequality with the conjugate pairs $(a_2,a_2')$ and $(a_3,a_3')$ we write
\begin{align*}
    & u(x)^{p+\sigma - 1} \le \frac{u(x)^{p^*_{s_1}+\sigma - 1}}{a_2} + \frac{1}{a_2'} \le u(x)^{p^*_{s_1} +\sigma - 1} +1, \text{ and } \\
    & u(x)^{q+\sigma - 1} \le \frac{u(x)^{p^*_{s_1}+\sigma - 1}}{a_3} + \frac{1}{a_3'} \le u(x)^{p^*_{s_1} +\sigma - 1} +1.
\end{align*}
Hence the R.H.S of \eqref{M to infty} can be estimated as
\begin{align*}
    \al\intom u(x)^{p+\sigma-1}\,\dx+\be\intom u(x)^{q+\sigma-1}\,\dx \le 2(\al + \be)(1+|\Om|) \left(1+ \intom u(x)^{p^*_{s_1} +\sigma - 1} \, \dx \right).
\end{align*}
Now using the facts $\sigma \ge 1$ and $\sigma +p -1 \le \sigma p$, we obtain from \eqref{M to infty} that 
\begin{align*}
    \left(1+ \int_{\rd}\left(u(x)^{\frac{\sigma+p-1}{p}}\right)^{p^*_{s_1}}\,\dx\right)^{\frac{p}{p^*_{s_1}}} \le C \left( \frac{\sigma + p -1}{p}\right)^{p-1} \left(1+ \intom u(x)^{p^*_{s_1} +\sigma - 1} \, \dx \right),
\end{align*}
where $C= C(\al, \be, \Om, d,s_1,p)>0$. Set $\vartheta=\sigma+p-1$. Then the above inequality can be written as 
\begin{align}\label{form 1}
    \left(1+ \int_{\rd} u(x)^{\frac{\vartheta}{p}p^*_{s_1}}\,\dx\right)^{\frac{p}{p^*_{s_1}(\vartheta - p)}} \le C^{\frac{1}{\vartheta - p}} \vartheta^{\frac{p-1}{\vartheta - p}} \left(1+ \int_{\rd} u(x)^{p^*_{s_1}+ \vartheta - p} \, \dx \right)^{\frac{1}{\vartheta-p}}.
\end{align}
We consider the sequences $(\vartheta_j)$ defined as follows
\begin{align*}
    \vartheta_1 = p^*_{s_1}, \vartheta_2 = p + \frac{p^*_{s_1}}{p}(\vartheta_1 -p), \cdot \cdot \cdot,  \vartheta_{j+1} = p + \frac{p^*_{s_1}}{p}(\vartheta_j -p).   
\end{align*}
Observe that $p^*_{s_1}-p+ \vartheta_{j+1} = \frac{p^*_{s_1}}{p} \vartheta_{j}$, and $\vartheta_{j+1} = p + \left(\frac{p^*_{s_1}}{p}\right)^j(\vartheta_1 - p)$. Since $p^*_{s_1} >p$, we get $\vartheta_{j} \ra \infty$, as $j \ra \infty$. From \eqref{form 1}, we then write 
\begin{align}\label{form 2}
    \left(1+ \int_{\rd} u(x)^{\frac{\vartheta_{j+1}}{p}p^*_{s_1}}\,\dx\right)^{\frac{p}{p^*_{s_1}(\vartheta_{j+1} - p)}} \le C^{\frac{1}{\vartheta_{j+1} - p}} \vartheta_{j+1}^{\frac{p-1}{\vartheta_{j+1} - p}} \left(1+ \int_{\rd} u(x)^{\frac{p^*_{s_1}}{p}\vartheta_{j}} \, \dx \right)^{\frac{p}{p^*_{s_1}(\vartheta_{j}-p)}}.
\end{align}
Set $D_j := \left(1+ \int_{\rd} u(x)^{\frac{p^*_{s_1}}{p}\vartheta_{j}} \, \dx \right)^{\frac{p}{p^*_{s_1}(\vartheta_{j}-p)}}.$ We iterate \eqref{form 2} to get 
\begin{align}\label{form 3}
    D_{j+1} \le \displaystyle C^{\sum_{k=2}^{j+1} \frac{1}{\vartheta_k - p}} \left( \prod_{k=2}^{j+1} \vartheta_{k}^{\frac{1}{\vartheta_{k} - p}}  \right)^{p-1} D_1,
\end{align}
where $D_1 = \left(1+ \int_{\rd} u(x)^{\frac{p^*_{s_1}}{p}p^*_{s_1}} \, \dx \right)^{\frac{p}{p^*_{s_1}(p^*_{s_1}-p)}}$ which is finite by using the claim, and 
\begin{align}\label{form 4}
    D_{j+1} \ge  \left( \left(\int_{\rd} u(x)^{\frac{p^*_{s_1} \vartheta_{j+1}}{p}} \, \dx \right)^{\frac{p}{p^*_{s_1}\vartheta_{j+1}}} \right)^{\frac{\vartheta_{j+1}}{\vartheta_{j+1}-p}} = \norm{u}_{L^{\frac{p^*_{s_1} \vartheta_{j+1}}{p}}(\rd)}^{\frac{\vartheta_{j+1}}{\vartheta_{j+1}-p}}.
\end{align} 
Combining \eqref{form 3} and \eqref{form 4} we have
\begin{equation}\label{form 5}
\norm{u}_{L^{\frac{p^*_{s_1} \vartheta_{j+1}}{p}}(\rd)}^{\frac{\vartheta_{j+1}}{\vartheta_{j+1}-p}}\leq C^{\sum_{k=2}^{j+1} \frac{1}{\vartheta_k - p}} \left( \prod_{k=2}^{j+1} \vartheta_{k}^{\frac{1}{\vartheta_{k} - p}}  \right)^{p-1} D_1.
\end{equation}
Moreover,
\begin{align*}
  & \sum_{k=2}^\infty\frac{1}{\vartheta_k-p}= \frac{1}{(\vartheta_1-p)} \sum_{k=2}^\infty \left(\frac{p}{p^*_{s_1}} \right)^{k-1} = \frac{p}{(p_{s_1}^*-1)(p_{s_1}^*-p)}, \text{ and } \\
  & \prod_{k=2}^{\infty}\vartheta_{k}^{\frac{1}{\vartheta_k - p}} = \exp\left( \sum_{k=2}^\infty \frac{\log(\vartheta_{k})}{\vartheta_{k} - p} \right) = \exp \left(\frac{p}{(p^*_{s_1} - p)^2} \log \left(p \left( \frac{p^*_{s_1}(p^*_{s_1}-p)}{p}  \right)^{p^*_{s_1}} \right) \right).
\end{align*}
Therefore, taking the limit as $j \ra \infty$ in \eqref{form 5}, we conclude that $u \in L^{\infty}(\rd)$.

\smallskip
\noi $\bm{d=s_1p\,}$\textbf{:} We proceed similarly as in the previous case by replacing the following fractional Sobolev inequality (whenever required):
\begin{equation*}
\iint\limits_{\rd\times\rd}\left|u_M(x)^{\frac{\sigma+p-1}{p}}-u_M(y)^{\frac{\sigma+p-1}{p}}\right|^{p}\,d\mu_1
\geq\Theta_{s_1,p}(\Omega)
\left(\int_{\rd}\left(u_M(x)^{\frac{\sigma+p-1}{p}}\right)^{2 p}\,\dx\right)^{\frac{1}{2}},
\end{equation*}
where 
$$
\Theta_{s_1,p}(\Omega):=\min_{\substack{u\in\wpso}}\left\{[u]^p_{s_1,p}:\norm{u}_{L^{2p}(\Omega)}=1\right\}.
$$
Following similar arguments as given in the case $d>s_1p$, we infer
\begin{align}\label{form 1.1}
    \left(1+ \int_{\rd} u(x)^{2 \vartheta}\,\dx\right)^{\frac{1}{2(\vartheta - p)}} \le C^{\frac{1}{\vartheta - p}} \vartheta^{\frac{p-1}{\vartheta - p}} \left(1+ \int_{\rd} u(x)^{p+ \vartheta} \, \dx \right)^{\frac{1}{\vartheta-p}}.
\end{align}
Then by considering the following sequences $(\vartheta_j)$ defined as:
\begin{align*}
    \vartheta_1 = 2p, \vartheta_2 = p + 2(\vartheta_1 -p), \cdot \cdot \cdot,  \vartheta_{j+1} = p + 2(\vartheta_j -p),   
\end{align*}
we obtain $u \in L^{\infty}(\rd)$. 
\\
\smallskip
\noi $\bm{d<s_1p\,}$\textbf{:} By the fractional Morrey's inequality (\cite[Proposition~2.9]{BrLiPa}), we see that functions in $\wpso$ are H\"older continuous and hence bounded. This completes the proof.
\end{proof}

We use the following version of the strong maximum principle for the positive solution of \eqref{main problem}. 

\begin{proposition}[Strong Maximum Principle]\label{SMP}
Let $\Om \subset \rd$ be a bounded open set and $0<s_2<s_1<1<q<p< \infty$. Let $u \in \wpso\cap L^\infty(\rd)$ be a non-negative supersolution of \eqref{main problem}. Then either $u>0$ a.e. in $\Om$ or $u \equiv 0$ a.e. in $\rd$.
\end{proposition}

\begin{proof} $\bm{\al, \be \ge 0}$\textbf{:} Since $u$ is a non-negative supersolution of \eqref{main problem}, we obtain 
\begin{align*}
 \left<A_p(u), v \right> + \left<B_q(u), v \right> \ge \al \intom u^{p-1}v + \be \intom  u^{q-1}v \ge 0,
\end{align*}
for every $v \in \wpso$ with $v \ge 0$. Now we can use \cite[(2) of Theorem 1.1]{Am} (by taking $c(x)=0$) with modifications (due to the presence of multiple parameters $s_1,s_2$) to conclude either $u>0$ a.e. in $\Om$ or $u \equiv 0$ a.e. in $\rd$.

\noi $\bm{\al, \be \le 0}$ or $\bm{\al\be\leq0}$\textbf{:} Let $x_0 \in \Om$ and $R>0$ be such that $B_R(x_0) \subset \Om$. Since $u$ is a non-negative supersolution of \eqref{main problem}, then we proceed as in \cite[Lemma 2.1]{Am}, for any $R_1 > 0$ satisfying $B_{R_1} = B_{R_1}(x_0) \subset B_{\frac{R}{2}}(x_0)$, and obtain the following logarithmic estimate 

\begin{equation}\label{logest1}
    \begin{split}
        \iint\limits_{B_{R_1} \times  B_{R_1}}\left| \log\left(\frac{u(x)+\delta}{u(y)+\delta}\right)\right|^q\, {\rm d}\mu_2 \le C \bigg(\delta^{1-q} R_1^d \bigg[ R^{-s_1 p}\,\Tail_{p}(u_-;x_0,R)^{p-1} \\ + R^{-s_2 q}\,\Tail_{q}(u_-;x_0,R)^{q-1} \bigg] 
        R_1^{d-s_1 p} + R_1^{d-s_2 q} \left(\norm{u}_{L^\infty(\rd)}+\delta\right)^{p-q}
        \\
        + \left(|\alpha|+|\beta|\,\norm{u}_{L^\infty(\rd)}^{p-q}\right)|B_{ 2 R_1 }(x_0)| \bigg),
    \end{split}
\end{equation}
where $\delta\in (0,1)$ and $C= C(d,s_1,p,s_2,q)>0$. Now the result follows using \eqref{logest1} and the arguments given in \cite[Lemma 2.3]{Am}.
\end{proof}

\section{Variational framework}\label{section3}

To obtain the existence part of Theorem \ref{Existence and nonexistence 1}-\ref{Existence main 2}, in this section, we study several properties of energy functionals associated with \eqref{main problem}. In view of Remark \ref{range changes}, we assume $s_2 < s_1$ and $q < p$ in the rest of the paper. We consider the following functional on $\wpso$:
\begin{align*}
    I_+(u) = \frac{\normo^p}{p} + \frac{\normt^q}{q} - \al\frac{\norm{u^{+}}^p_p }{p} - \be \frac{\norm{u^{+}}^q_q}{q}, \quad \forall \, u \in \wpso.
\end{align*}
Now we define
\begin{align*}
    &\left<A_p(u), \phi\right> = \iint \limits_{\rd \times \rd} |u(x)-u(y)|^{p-2}(u(x)-u(y))(\phi(x)-\phi(y)) \, \dmuo; \\
    &\left<B_q(u), \phi\right> = \iint \limits_{\rd \times \rd} |u(x)-u(y)|^{q-2}(u(x)-u(y))(\phi(x)-\phi(y)) \, \dmut, \quad \forall \, u,\phi \in \wpso,
\end{align*}
where $\left< \cdot \right>$ denotes the duality action. Using the H\"{o}lder's inequality, it follows that $\norm{A_p(u)} \le [u]_{s_1, p}^{p-1}$ and $\norm{B_q(u)} \le [u]_{s_2, q}^{q-1}$. One can verify that $I_+ \in C^1(\wpso, \R)$ and
\begin{align*}
    \left<I_+'(u), \phi \right> = \left<A_p(u), \phi\right> + \left< B_q(u), \phi \right> - \al \intom (u^{+})^{p-1} \phi \,\dx - \be \intom (u^{+})^{q-1} \phi \, \dx, \quad \forall \, u, \phi \in \wpso.
\end{align*}
\begin{remark}\label{nonnegative}
If $u \in \wpso$ is a critical point of $I_+$, i.e., $\big<I_+'(u), \phi \big> =0$ for all $\phi \in \wpso$, then $u$ is a solution of \eqref{main problem}. Moreover, for $\phi=-u^{-}$, using (i) of Lemma \ref{inequality} we see 
\begin{align*}
   0=\left<I_+'(u), -u^{-} \right> = \left<A_p(u), -u^{-} \right> + \left<B_q(u), -u^{-} \right> \ge [u^{-}]^p_{s_1,p} + [u^{-}]^q_{s_2,q}.
\end{align*}
The above inequality yields $u^{-} = c$ a.e. in $\rd$ for some $c\in\re$. Moreover, since $u^{-}\in\wpso$, we get $c=0$.  Thus every critical point of $I_+$ is a nonnegative solution of \eqref{main problem}.
\end{remark}

Now we discuss the coercivity and weak lower semicontinuity of $I_+$. 

\begin{proposition}\label{coercive0}
Let $\al < \lao$ and $\be >0$. Then the functional $I_+$ is weakly sequentially lower semicontinuous, coercive, and bounded below on $\wpso$.
\end{proposition}

\begin{proof}
Let $u_n \rightharpoonup u$ in $\wpso$. Then using the compactness of the embeddings $\wpso \hookrightarrow L^p(\Om)$, $\wqst \hookrightarrow L^{q}(\Om)$, and the weak lower semicontinuity of the seminorm, we get 
\begin{align*}
    \lowlim_{n \ra \infty} {I_+(u_n)}  = \lowlim_{n \ra \infty} \frac{[u_n]_{s_1,p}^p}{p} + \lowlim_{n \ra \infty} \frac{[u_n]_{s_2,q}^q}{q} - \al \lim_{n \ra \infty} \frac{\norm{u_n^+}^p_p }{p} - \be \lim_{n \ra \infty}  \frac{\norm{u_n^+}^q_q}{q} \ge I_+(u).
\end{align*}
Now we prove the coercivity of $I_+$. Suppose $\al \le 0$. Then using $\wpso \hookrightarrow L^q(\Om)$, 
\begin{align}\label{co1}
    I_+(u) \ge \frac{[u]_{s_1,p}^p}{p} - \be \frac{\norm{u^+}^q_q}{q} \ge \frac{[u]_{s_1,p}^p}{p} - C \be \frac{[u]_{s_1,p}^q}{q}, \quad \forall \, u \in \wpso \setminus \{ 0\}.
\end{align}
If $\al > 0$, then there exists $a \in (0,1)$ such that $\al = a \lao$. In this case, using  $\wpso \hookrightarrow L^q(\Om)$, we get 
\begin{align}\label{co1.1}
    I_+(u) \ge \frac{[u]_{s_1,p}^p}{p} - a \lao \frac{\norm{u^+}_p^p}{p} - \be \frac{\norm{u^+}^q_q}{q} \ge \frac{[u]_{s_1,p}^p}{p} - a \lao \frac{\norm{u^+}_p^p}{p} - C \be \frac{[u]_{s_1,p}^q}{q},
\end{align}
for every $u \in \wpso \setminus \{0\}$. From the definition of $\lao$, we have $[u]_{s_1,p}^p \ge \lao \norm{u}_p^p \ge \lao \norm{u^+}_p^p$. Therefore, \eqref{co1} yields
\begin{align}\label{co2}
    I_+(u) \ge \frac{1-a}{p} [u]_{s_1,p}^p - C \be \frac{[u]_{s_1,p}^q}{q}, \quad \forall \, u \in \wpso \setminus \{ 0 \}.
\end{align}
In view of \eqref{co1}, observe that \eqref{co2} holds for every $\al < \lao$. For any $\ep >0$, applying Young's inequality with the conjugate pair $(\frac{p}{q}, \frac{p}{p-q})$ we obtain 
\begin{align*}
    [u]_{s_1,p}^q \le \ep \frac{q}{p} [u]_{s_1,p}^p + \frac{p-q}{p} \ep^{-\frac{p}{p-q}}. 
\end{align*}
Hence from \eqref{co1} we have the following estimate for every $u \in \wpso \setminus \{ 0\}$: 
\begin{align*}
    I_+(u) \ge  \frac{1-a}{p}  [u]_{s_1,p}^p  - \ep \frac{C \be}{p} [u]_{s_1 , p}^p - \frac{C \be(p-q)}{qp} \ep^{-\frac{p}{p-q}}. 
\end{align*}
We choose $\ep>0$ so that $C \be \ep < \frac{1-a}{2}$. Therefore, from the above estimate, we get 
\begin{align*}
    I_+(u) \ge \frac{1-a}{2p} [u]_{s_1,p}^p - \frac{C \be(p-q)}{qp} \ep^{-\frac{p}{p-q}} \quad \forall \, u \in \wpso \setminus \{ 0\}.
\end{align*}
Thus the functional $I_+$ is coercive on $\wpso$. Next, we prove that $I_+$ is bounded below. Set $M>0$ such that $M^{p-q} \ge p(1+C \be q^{-1})$. Then using \eqref{co1}, we get 
\begin{align*}
    I_+(u) \ge [u]_{s_1,p}^q \left( \frac{[u]_{s_1,p}^{p-q}}{p} -  \frac{C \be}{q} \right) \ge M^q, \; \text{ provided } \; [u]_{s_1,p} \ge M. 
\end{align*}
Further, if $[u]_{s_1,p} \le M$, then $I_+(u) \ge - \frac{ M^q C \be}{q}$. Thus, $I_+$ is bounded below on $\wpso$.   
\end{proof}

In the following proposition, we verify that $I_+$ satisfies the Palais-Smale (P.S.) condition on $\wpso$. 

\begin{proposition}\label{PS}
Let $\al \neq \la^1_{s_1,p}$. Let $(u_n)$ be a sequence in $\wpso$ such that $I_+(u_n) \ra c$ for some $c \in \R$ and $I_+'(u_n) \ra 0$ in $(\wpso)^*$. Then $(u_n)$ possesses a convergent subsequence in $\wpso$. 
\end{proposition}

\begin{proof}
First, we show that the sequence $(u_n)$ is bounded in $\wpso$. On a contrary, assume that $[u_n]_{s_1,p} \ra \infty$, as $n \ra \infty$. Using (i) of Lemma \ref{inequality}, note that 
\begin{align*}
 [u_n^-]^p_{s_1,p} \le [u_n^-]^p_{s_1,p} + [u_n^-]^q_{s_2,q} \le  \left| \left<I_+'(u_n), -u_n^{-} \right>\right| \le \norm{I_+'(u_n)} [u_n^-]_{s_1,p}.
\end{align*}
Hence $[u_n^-]_{s_1, p} \ra 0$, as $n \ra \infty$. Set $w_n=u_n [u_n]_{s_1,p}^{-1}$. Up to a subsequence, $w_n \rightharpoonup w$ in $\wpso$ and by the compactness of $\wpso \hookrightarrow L^p(\Om)$, $w_n \ra w$ in $L^p(\Om)$. Further, $[w_n^-]_{s_1,p} = [u_n^-]_{s_1,p} [u_n]_{s_1,p}^{-1} \ra 0$, as $n \ra \infty$. Therefore, $w_n^- \ra 0$ in $\wpso$ and hence in $L^p(\Om)$. This implies that $w_n^+ \ra w$ in $L^p(\Om)$, which yields $w \ge 0$ a.e. in $\Om$. We show that $w$ is an eigenfunction of the fractional $p$-Laplacian corresponding to $\al$. For any $\phi \in \wpso$, we write
\begin{align}\label{PS0}
    \left<A_p(u_n), \phi\right> + \left< B_q(u_n), \phi \right> - \al \intom \abs{u_n}^{p-2}u_n \phi - \be \intom \abs{u_n}^{q-2}u_n \phi =\ep_n, 
\end{align}
where $\ep_n \ra 0$ as $n \ra \infty$. From the above inequality, we obtain 
\begin{align}\label{PS1}
    \left<A_p(w_n),\phi\right> + [u_n]_{s_1,p}^{q-p} \left< B_q(w_n), \phi \right> - \al \intom \abs{w_n}^{p-2} w_n \phi - \be [u_n]_{s_1,p}^{q-p} \intom \abs{w_n}^{q-2} w_n \phi = \frac{\ep_n}{[u_n]_{s_1,p}^{p-1}}.  
\end{align}
Using the H\"{o}lder's inequality with the conjugate pair $(q, q')$, the Poincar\'{e} inequality $\norm{\phi}_q \le C(\Omega ) [\phi]_{s_1,p}$, and the boundedness of $(w_n)$ in $\wqst$ we have
\begin{align*}
    \abs{\left< B_q(w_n), \phi \right>} \le [w_n]_{s_2,q}^{q-1} [\phi]_{s_2,q} \le C [\phi]_{s_1,p}, \text{ and } \intom \abs{w_n}^{q-1} \abs{\phi} \le \norm{w_n}_{q}^{q-1}\norm{\phi}_{q} \le C[\phi]_{s_1,p},
\end{align*}
We choose $\phi=w_n-w$ in \eqref{PS1}, and take the limit as $n \ra \infty$ to get $\left<A_p(w_n), w_n-w\right> \ra 0.$ Further, since $A_p$ is a continuous functional on $\wpso$, we also have $ \left<A_p(w), w_n-w\right> \ra 0$. Further, using the definition of $A_p$
\begin{align}\label{PS2}
    \left<A_p(w_n) - A_p(w), w_n-w\right> \ge \left( [w_n]^{p-1}_{s_1,p} - [w]^{p-1}_{s_1,p} \right) \left( [w_n]_{s_1,p} - [w]_{s_1,p} \right).
\end{align}
Therefore, $[w_n]_{s_1,p} \ra [w]_{s_1,p}$, and hence the uniform convexity of $\wpso$ ensures that $w_n \ra w$ in $\wpso$. Further, since $[w]_{s_1,p}=1$ we also have $w \neq 0$ in $\Om$. Now using \eqref{PS1}, we obtain 
\begin{equation*}
    \left< A_p(w), \phi \right> = \al \intom |w|^{p-2}w \phi, \quad \forall \, \phi \in \wpso.
\end{equation*}
Thus $w$ is a nonnegative  weak solution to the problem 
\begin{equation}
\begin{aligned} 
 (-\De)^{s_1}_pu=\alpha|u|^{p-2}u \;\;\text{in }\Omega, \; u=0 \;\;\text{in } \Om^c.
 \end{aligned}
 \end{equation}
Now by the strong maximum principle for fractional $p$-Laplacian \cite[Proposition 2.6]{BrPa}, we conclude that $w>0$ a.e. in $\Om$. Therefore, the uniqueness of $\la^1_{s_1,p}$ (Proposition \ref{first eigenvalue}) yields $\al = \la^1_{s_1,p}$, resulting in a contradiction. Thus, the sequence $(u_n)$ is bounded in $\wpso$. By the reflexivity, up to a subsequence, $u_n \rightharpoonup \tilde{u}$ in $\wpso$. By taking $\phi=u_n-\tilde{u}$ in \eqref{PS0} and using the compact embeddings of $\wpso \hookrightarrow L^{\ga}(\Om)$ with $\ga \in [1,p]$, we get $\left<A_p(u_n), u_n-\tilde{u} \right> + \left<B_q(u_n), u_n - \tilde{u} \right> \ra 0$. Therefore, $\left<A_p(u_n)- A_p(\tilde{u}), u_n-\tilde{u} \right> + \left<B_q(u_n)- B_q(\tilde{u}), u_n - \tilde{u} \right> \ra 0$, which implies $\left<A_p(u_n)- A_p(\tilde{u}), u_n-\tilde{u} \right>\to 0$. Thus, $[u_n]_{s_1,p} \ra [\tilde{u}]_{s_1,p}$ (by using \eqref{PS2}), and from the uniform convexity, $u_n \ra \tilde{u}$ in $\wpso$, as required.
\end{proof}
The following lemma discusses the mountain pass geometry of $I_+$ for certain ranges of $\al$ and $\be$.

\begin{lemma}\label{MPG}
Let $\al> \la^1_{s_1,p}$ and $\be \le \al$. For $\rho>0$, let
\begin{align*}
    \S_{\rho}= \left\{u: \wpso : [u]_{s_1,p} = \rho \right\}.
\end{align*}
The following hold: 
\begin{enumerate}
    \item[(i)] There exist $\de= \de(\rho)>0$, and $\al_1 =\al_1(\rho)>0$ such that if $\al \in (0, \al_1)$, then $I_+(u) \ge \de$ for every $u\in \S_{\rho}$.
    \item[(ii)] There exists $v \in \wpso$ with $[v]_{s_1,p} > \rho$ such that $I_+(v)<0$.
\end{enumerate}
\end{lemma}

\begin{proof}
\noi (i) Let $\rho>0$ and $u \in \S_{\rho}$. Then using $\wpso \hookrightarrow L^\ga(\Om)$ for $\ga \in [1,p]$,  
\begin{equation}\label{MPG1}
    I_+(u) \ge \frac{\normo^p}{p} - \al \frac{\norm{u^+}^p_p}{p}  - \be \frac{\norm{u^+}^q_q}{q}  \ge \normo^q \left( \frac{\normo^{p-q}}{p} - C \al \frac{\normo^{p-q}}{p}  - C\frac{\al}{q}\right) = \rho^{q} A(\rho), 
\end{equation}
where $A(\rho)=\frac{\rho^{p-q}}{p} - C \al \frac{\rho^{p-q}}{p}  - C\frac{\al}{q}$. Choose $ 0 < \al_1 < \frac{\rho^{p-q}}{p} \left( C\frac{\rho^{p-q}}{p} + \frac{C}{q} \right)^{-1}$ and $\de = \rho^{q} A(\rho)$ with $\al \in (0, \al_1)$. Therefore, from \eqref{MPG1}, $I_+(u) \ge \de$ for all $\al \in (0, \al_1)$.

\noi (ii) Note that  
\begin{equation*}
    I_+(t\phi_{s_1,p}) = \frac{t^p}{p} \left([\phi_{s_1,p}]_{s_1,p}^p - \al \norm{\phi_{s_1,p}}_p^p \right) + \frac{t^q}{q} \left([\phi_{s_1,p}]_{s_2,q}^q - \be \norm{\phi_{s_1,p}}_q^q \right).
\end{equation*}
Since $p > q$ and $\al > \la^1_{s_1,p}$, we obtain $I_+(t\phi_{s_1,p}) \ra -\infty$, as $t \ra \infty$. Hence there exists  $t_1 > \rho [\phi_{s_1,p}]_{s_1,p}^{-1}$ such that $I_+(t\phi_{s_1,p})<0$ for all $t \ge t_1$. Thus, $v=t\phi_{s_1,p}$ with $t>t_1$ is the required function.
\end{proof}

\subsection{Nehari manifold}
This subsection briefly discusses the Nehari manifold associated with \eqref{main problem} and some of its properties. 

\begin{definition}[Nehari Manifold]
We define the Nehari manifold associated with \eqref{main problem} as
\begin{align*}
    \mathcal{N}_{\al,\be} := \left\{u \in \wpso \setminus \{0\} : \left<I_+'(u), u\right> =0\right\}.
\end{align*}
\end{definition}
Note that every nonnegative solution of \eqref{main problem} lies in $\mathcal{N}_{\al, \be}$. Now we provide a sufficient condition for which every critical point in $\N_{\al, \be}$ becomes a nonnegative solution of \eqref{main problem}. We consider the following functionals on $\wpso$:
\begin{align*}
    H_{\al}(u) = \normo^p - \al \norm{u^+}^p_p, \text{ and }  G_{\be}(u) = \normt^q - \be \norm{u^+}^q_q, \quad \forall \, u \in \wpso.
\end{align*}
Clearly, $H_{\al}, G_{\be} \in C^1(\wpso, \R)$, and the identity $\left<I_+'(u), u\right> = H_{\al}(u) + G_{\be}(u)$ holds.  
\begin{proposition}\label{Nehari0}
Assume that either $H_{\al}(u) \neq 0$ or $G_{\be}(u) \neq 0$. If $u$ is a critical point of $I_+$ in $\N_{\al, \be}$, then $u$ is a critical point of $I_+$ in $\wpso$. 
\end{proposition}
\begin{proof}
The proof follows using the arguments given in \cite[Lemma~2]{BoTa1}.
\end{proof}

Next, we state a condition for the existence of a critical point in $\N_{\al, \be}$. Let $H_{\al}(u),G_{\be}(u) \neq 0$ for some $u \in \wpso$. Define 
$$ t_{\al,\be} (= t_{\al,\be}(u)) := \left(-\frac{G_{\be}(u)}{H_{\al}(u)} \right)^{\frac{1}{p-q}}.$$ 

Notice that, for $t \in \R$, $\left<I_+'(t u), tu \right> = t (t^{p-1}H_{\al}(u) +  t^{q-1}G_{\be}(u))$. In particular, $t_{\al,\be}u \in \N_{\al, \be}$.  

\begin{proposition}\label{Nehari1}
Let $u \in \wpso$. The following hold:
\begin{enumerate}
    \item[(i)] If $G_{\be}(u)<0<H_{\al}(u)$, then $I_+(t_{\al,\be}u)= \underset{t \in \R^+}{\min} I_+(tu)$, and $I_+(t_{\al,\be}u)<0$. Moreover, $t_{\al,\be}$ is unique.
    \item[(ii)] If $H_{\al}(u) < 0 < G_{\be}(u)$, then $I_+(t_{\al,\be}u)= \underset{t \in \R^+}{\max} \, I_+(tu)$, and $I_+(t_{\al,\be}u)>0$. Moreover, $t_{\al,\be}$ is unique.
\end{enumerate}
\end{proposition}

\begin{proof}
The proof follows using the same arguments presented in \cite[Proposition~6]{BoTa1}.
\end{proof}

\begin{remark}\label{Halpha1}
(i) Let $u \in \N_{\al, \be}$. Then $H_{\al}(u) + G_{\be}(u) =0$, and hence 
\begin{align*}
    I_+(u)= \frac{p-q}{pq}G_{\be}(u)=  \frac{q-p}{pq} H_{\al}(u).
\end{align*}
From the above identity, it is clear that if $I_+(u) \neq 0$, then either $G_{\be}(u)<0<H_{\al}(u)$ or $H_{\al}(u) < 0 < G_{\be}(u)$. 

\noi (ii) If $u \in \N_{\al,\be}$, and $H_{\al}, G_{\be}$ satisfy the assumptions given in the above proposition, then from (i) and Proposition \ref{Nehari1}, $t_{\al,\be} = 1$ is the unique minimum or maximum point on $\R^+$. 
\end{remark}

\begin{remark}\label{Halpha}
Using Proposition \ref{first eigenvalue} and $\norm{u^+}_{\ga} \le \norm{u}_{\ga}$, we get
\begin{enumerate}
    \item[(i)] if $\al < \lao$, then $H_{\al}(u) > \normo^p - \lao \norm{u^+}^p_p \ge \normo^p - \lao \norm{u}^p_p \ge 0$ for $u \in \wpso \setminus \{0\}$,
    \item[(ii)] if $\be < \lat$, then $G_{\be}(u) > \normo^p - \lat \norm{u^+}^q_q \ge \normt^q - \lat \norm{u}^q_q \ge 0$ for $u \in \wpso \setminus \{0\}$.
\end{enumerate}
\end{remark}
\subsection{Method of sub and super solutions}
We consider the following energy functional on $\wpso$:
\begin{align*}
     I(u) = \frac{\normo^p}{p} + \frac{\normt^q}{q} - \al\frac{\norm{u}^p_p }{p} - \be \frac{\norm{u}^q_q}{q}, \quad \forall \, u \in \wpso.
\end{align*}
Notice that $I \in C^1(\wpso, \R)$, and  
\begin{align*}
    \left<I'(u), \phi\right> = \left<A_p(u), \phi\right> + \left< B_q(u), \phi \right> - \al \intom \abs{u}^{p-2}u \phi \,\dx - \be \intom \abs{u}^{q-2}u \phi, \quad \forall \, u, \phi \in \wpso.
\end{align*}
In this subsection, using sub and super solutions techniques, we discuss the existence of critical points for $I$. We say $\overline{u} \in \wpso$ is a supersolution of \eqref{main problem}, if 
\begin{equation}\label{subsuper1}
    \left<A_p(\overline{u}), \phi \right>  + \left<B_q(\overline{u}) , \phi \right> \ge \al \intom \abs{\overline{u}}^{p-2}\overline{u}\phi \,\dx + \be \intom \abs{\overline{u}}^{q-2}\overline{u}\phi \, \dx, \quad \forall \, \phi \in \wpso,\, \phi \ge 0.
\end{equation}
A function $\underline{u} \in \wpso$ is called a subsolution of \eqref{main problem} if the reverse inequality holds in \eqref{subsuper1}. 

\begin{definition}[Truncation function]
    Let $\underline{u}, \overline{u} \in L^{\infty}(\Omega)$ be such that $\underline{u} \le \overline{u}$ a.e. in $\Om$. For $t \in \R$, we define the truncation function corresponds to $f(t)= \al \abs{t}^{p-2}t + \be \abs{t}^{q-2}t$ as follows:
\begin{equation}\label{truncation}
 \tilde{f}(x,t) := \left\{\begin{array}{ll}
                 f(\overline{u}(x)) & \text{ if } t \ge \overline{u}(x),\\
                 f(t) & \text{ if } \underline{u}(x) \le t \le \overline{u}(x), \\
                 f(\underline{u}(x)) & \text{ if } t \le \underline{u}(x).
          \end{array}\right.
\end{equation}
\end{definition}
By definition, $\tilde{f}(\cdot, t)$ is continuous on $\R$. Further, using $\underline{u}, \overline{u} \in L^{\infty}(\Omega)$ it is easy to see that $\tilde{f} \in L^{\infty}(\Omega \times \R)$. Now we consider the following functional associated with $\tilde{f}(\cdot, u(x))$: 
\begin{equation*}
    \tilde{I}(u) =  \frac{\normo^p}{p} + \frac{\normt^q}{q} - \intom \tilde{F}(x,u(x)) \, \dx, \quad \forall \, u \in \wpso,  
\end{equation*}
where $\tilde{F}( x, u(x)) := \int_0^{u(x)} \tilde{f} (x, \tau) \, \dtau$. 
Note that, for $u(x) \in (\underline{u}(x), \overline{u}(x))$, $\tilde{I}$ coincides with the energy functional $I$. Further, $\tilde{I} \in \C^1(\wpso, \R)$, and 
\begin{equation*}
    \left< (\tilde{I})'(u), \phi \right> = \left<A_p(u), \phi\right> + \left< B_q(u), \phi \right> - \intom \tilde{f}(x,u(x)) \phi(x) \, \dx, \quad \forall \, u, \phi \in \wpso. 
\end{equation*}
In the following proposition, we prove some properties of $\tilde{I}$ that ensure the existence of critical points for $\tilde{I}$.
\begin{proposition}\label{critical point of truncation}
Let $\underline{u}, \overline{u} \in L^{\infty}(\Om)$ be such that $\underline{u} \le \overline{u}$ a.e. on $\Om$. Then $\tilde{I}$ is bounded below, coercive and weak lower semicontinuous on $\wpso$.
\end{proposition}
\begin{proof}
Since $\underline{u}, \overline{u} \in L^{\infty}(\Om)$, there exists $C>0$ such that 
$\abs{\tilde{f}(x,t)} \le C$, and $\abs{\tilde{F}(x,t)} \le C|t|$, for all $x\in \Om, t \in \R$. Hence for $u \in \wpso$
\begin{align*}
    \tilde{I}(u) \ge \frac{\normo^p}{p} + \frac{\normt^q}{q} - C \norm{u}_1 \ge  \frac{\normo^p}{p} + \frac{\normt^q}{q} - C [u]_{s_2,q} \abs{\Om}^{\frac{1}{q'}}.
\end{align*}
Now using similar arguments as in Proposition \ref{coercive0}, it follows that $\tilde{I}$ is coercive and bounded below on $\wpso$. Next, for a sequence $u_n \rightharpoonup u$ in $\wpso$,
\begin{equation}\label{semi1}
    \lowlim_{n \ra \infty} \tilde{I}(u_n) \ge \frac{\normo^p}{p} + \frac{\normt^q}{q} - \lim_{n \ra \infty} \intom \tilde{F}(x,u_n(x)) \, \dx.
\end{equation}
We claim that $\int_{\Om} \tilde{F}(x,u_n(x)) \, \dx \ra \int_{\Om} \tilde{F}(x,u(x)) \, \dx$.
By the compact embeddings of $\wpso \hookrightarrow L^p(\Omega )$, we have $u_n \ra u$ in $L^p(\Om)$ and hence $u_n \ra u $ in $L^1(\Omega )$. 
Further, using $\tilde{f} \in L^{\infty}(\Omega \times \R )$, 
\begin{align*}
    \left| \int_{\Omega }  \left( \tilde{F}(x,u_n(x)) - \tilde{F}(x,u(x)) \right) \, \dx \right| \le  \int_{\Omega} \int_{u(x)}^{u_n(x)} \abs{\tilde{f}(x, \tau} \, \dtau \dx \le M \int_{\Om} \abs{u_n(x) - u(x)}  \, \dx, 
\end{align*}
and the claim follows. Therefore, in view of \eqref{semi1}, $\tilde{I}$ is weak lower semicontinuous on $\wpso$.
\end{proof}
In the following proposition, we prove that every critical point of $\tilde{I}$ lies between sub and super solutions. 
\begin{proposition}\label{sub and sup}
Let $\underline{u}, \overline{u} \in L^{\infty}(\Om)$ be such that $\underline{u} \le \overline{u}$ a.e. in $\rd$. If  $u \in \wpso$ is a critical point of $\tilde{I}$, then $\underline{u} \le u \le \overline{u}$ a.e. in $\rd$.
\end{proposition}

\begin{proof}
From the definition of sub and super solutions, it is clear that $\underline{u} = u = \overline{u} =0$ in $\rd \setminus \Om$,  since each function lies in $\wpso$. Now we show that $\underline{u} \le u \le \overline{u}$ a.e. in $\Omega$. Our proof is by the method of contradiction. On the contrary, assume that $u \ge \overline{u}$ on $A \subset \Omega$ with $|A| >0$. We choose $(u - \overline{u})^+ \in \wpso$ as a test function. Using $\overline{u}$ is a supersolution of \eqref{main problem} and $u$ is a critical point of $\tilde{I}$, together with \eqref{truncation} we get
\begin{align*}
    & \left< A_p(\overline{u}), (u - \overline{u})^+ \right>  + \left<B_q(\overline{u}) , (u - \overline{u})^+ \right> \ge \al \intom \abs{\overline{u}}^{p-2}\overline{u}(u-\overline{u}) + \be \intom \abs{\overline{u}}^{q-2}\overline{u}(u-\overline{u}),
    \\
    & \left< A_p(u), (u - \overline{u})^+ \right>  + \left<B_q(u) , (u - \overline{u})^+ \right> = \intom f(\overline{u}) (u-\overline{u}) = \intom\left(\al \abs{\overline{u}}^{p-2} + \be\abs{\overline{u}}^{q-2}\right)\overline{u}(u-\overline{u}).
\end{align*}
The above inequalities yield
\begin{align}\label{sub and sup 1}
    \big<A_p(u)-A_p(\overline{u}), (u - \overline{u})^+ \big> + \big<B_q(u)-B_q(\overline{u}), (u - \overline{u})^+ \big> \le 0.
\end{align} 
From the definition of $A_p$,
\begin{equation*}
    \begin{split}
        \left<A_p(u)-A_p(\overline{u}), (u - \overline{u})^+ \right> = \iint \limits_{\rd \times \rd} \left( |u(x)-u(y)|^{p-2}(u(x)-u(y)) - |\overline{u}(x)-\overline{u}(y)|^{p-2}(\overline{u}(x)-\overline{u}(y)) \right) \\
       \left((u(x) -\overline{u}(x))^+ - (u(y)- \overline{u}(y))^+ \right) \, \dmuo. 
       \end{split}
\end{equation*}
Now we consider the following cases:

\noi $\bm{2 \le q <p}$\textbf{:} Without loss of generality, we assume that $u(x) -\overline{u}(x) \ge u(y)- \overline{u}(y)$. Otherwise, exchange the roll of $x$ and $y$. 
Applying (ii) and (i) of Lemma \ref{inequality}, we then obtain 
\begin{equation*}
    \begin{split}
     \left<A_p(u)-A_p(\overline{u}), (u-\overline{u})^+ \right>
     & \ge C(p)  \iint \limits_{\rd \times \rd}  \abs{ (u(x) - \overline{u}(x)) -(u(y) -\overline{u}(y))}^{p-2}  \\ 
     & \quad \left( (u(x) -\overline{u}(x)) - (u(y)- \overline{u}(y) \right)
    \left((u(x) -\overline{u}(x))^+ - (u(y)- \overline{u}(y))^+ \right) \, \dmuo\\
    & \ge C(p) [(u-\overline{u})^+]^p_{s_1,p}.
    \end{split}
\end{equation*}
Similarly, we can show that $\big<B_q(u)-B_q(\overline{u}), (u-\overline{u})^+ \big> \ge C(q)[(u-\overline{u})^+]^q_{s_2,q}$. Therefore, from \eqref{sub and sup 1}, $[(u-\overline{u})^+]_{s_1,p} = 0$. By Poincar\`{e} inequality, $\norm{(u-\overline{u})^+}_p \le C[(u-\overline{u})^+]_{s_1,p}= 0$, which is a contradiction. 

\noi $\bm{q < 2 \le p}$\textbf{:} In this case, using \cite[Lemma 2.4]{IaMoSq1} (for $B_q$) we obtain, 
\begin{align*}
     & \left< A_p(u) - A_p(\overline{u}), (u - \overline{u})^+ \right> \ge C(p) [(u-\overline{u})^+]^p_{s_1,p} ; \\
     & \left< B_q(u) - B_q(\overline{u}), (u - \overline{u})^+ \right> \ge C(q) \frac{[(u-\overline{u})^+]^2_{s_2,q}}{\left([u]^q_{s_2,q} + [\overline{u}]^q_{s_2,q} \right)^{\frac{2-q}{q}}}.
\end{align*}
Hence, we get a contradiction using \eqref{sub and sup 1}. For $q < p < 2$, again using \cite[Lemma 2.4]{IaMoSq1}, we similarly get 
a contradiction. Thus  $u \le \overline{u}$ a.e. in $\rd$. Now suppose $u \le \underline{u}$ in $A \subset \Omega$ with $|A| > 0$, then taking $(u-\underline{u})^- \in \wpso$ as a test function, we also get a contradiction for all possible choices of $p$ and $q$. Therefore, $\underline{u} \le u \le \overline{u}$ a.e. in $\rd$. 
\end{proof}

\section{Existence and non-existence of positive solutions}\label{Section4}

Depending on the ranges of $\al, \be$, this section is devoted to proving the existence and non-existence of positive solutions for \eqref{main problem}. This section's terminology `solution' is meant to be nontrivial unless otherwise specified. First, we consider the region where $\al, \be$ do not exceed $\lao, \lat$ respectively.

\begin{proposition}\label{Nonexistence1} 
It holds
\begin{itemize}
    \item[(i)] Let $(\al,\be) \in \left((-\infty, \lao) \times (-\infty, \lat)\right) \cup \left( \{\lao\} \times (-\infty , \lat) \right) \cup \left( (-\infty, \lao) \times \{\lat\} \right)$. Then \eqref{main problem} does not admit a solution.
    \smallskip
    \item[(ii)] Let $\al = \lao$ and $\be = \lat$. Then \eqref{main problem} admits a solution if and only if \eqref{LI} violates.
\end{itemize}
\end{proposition}
\begin{proof}
(i) Let $\al < \lao$ and $\be < \lat$. Suppose  $u \in \wpso \setminus \{0\}$ is a solution of \eqref{main problem}. Then using the definition of $\lao$ and $\lat$ (Proposition \ref{first eigenvalue}), we get 
\begin{align}\label{nonexists0}
    0 < (\lao - \al)\norm{u}^p_p \le \normo - \al \norm{u}^p_p = \be \norm{u}^q_q - \normt \le (\be - \lat)\norm{u}^q_q < 0.  
\end{align}
A contradiction. Therefore, \eqref{main problem} does not admit a solution. For other cases, contradiction similarly follows using \eqref{nonexists0}.

\noi (ii) For $\al = \lao$ and $\be = \lat$, if $u \in \wpso \setminus \{0\}$ is a solution of \eqref{main problem}, then the equality occurs in \eqref{nonexists0}. As a consequence, $u$ becomes an eigenfunction corresponding to both $\lao$ and $\lat$, i.e., \eqref{LI} violates. Conversely, suppose \eqref{LI} does not hold. For $\al \le \lao$ and $\be \le \lat$, using Remark \ref{Halpha} we have $I_+(u) \ge 0$  for any $u \in \wpso \setminus \{0\}$. Thus $0$ is the global minimizing point for $I_+$. Further, since $\phi_{s_1,p} = c\phi_{s_2,q}$ for some nonzero $c \in \R$, by setting $\tilde{u}=c_1\phi_{s_1,p} = c_2\phi_{s_2,q}$ (where $c_1, c_2 \neq 0$) we see that $I_+(\tilde{u})=0$. Therefore, $\tilde{u} \neq 0$ is a solution of \eqref{main problem}.  
\end{proof}

Before going to the proof of theorem \ref{Existence and nonexistence 1}, we recall a result from \cite{NgVo}, where for $d>s_1p$ the authors provided the existence of a positive solution of \eqref{weight problem1}. However, we stress that the same conclusion can be drawn for $d\leq s_1p$. For $0<s<1\le r< \infty$ and $m_r \in L^{\infty}(\Om)$ with $m_r^+ \not \equiv 0$, we denote
\begin{align*}
    \lambda^1_{s,r}(\Om, m_r) := \inf \left\{ [u]_{s,r}^r : u \in \wrs \text{ and } \intom m_r|u|^r=1 \right\}
\end{align*}
as the first Dirichlet eigenvalue of the weighted eigenvalue problem of the fractional $r$-Laplace operator (see \cite{DeQu}).
\begin{theorem}[\protect{\cite[Theorem 1.1]{NgVo}}]\label{Nguyen}
Let $\Om \subset \rd$ be a bounded open set, $0<s_2< s_1< 1< q \le p< \infty$, and $m_p, m_q \in L^{\infty}(\Om)$ with $m_p^+, m_q^+ \not \equiv 0$. Let $\lao(\Om,m_p), \lat(\Om,m_q)$ be respectively the first Dirichlet eigenvalue of weighted eigenvalue problems for fractional $p$-Laplace and fractional $q$-Laplace operators with weights $m_p, m_q$. Suppose, $\lao(\Om, m_p) \neq \lat(\Om, m_q)$. Then for $$\al > \min \{\lao(\Om,m_p), \lat(\Om,m_q) \},$$ the problem \eqref{weight problem1} admits a positive solution.
\end{theorem}

\noi \textbf{Proof of Theorem \ref{Existence and nonexistence 1}:}
(i) $\bm{\al > \lao,\,\be < \lat\,}$\textbf{:} Let $\be > 0$. Then using $\al > \lao$ and $\be < \lat$, we get 
\begin{align*}
  \lao \left(\Om, \frac{\al}{\be} \right) =  \frac{\lao }{\al}\be < \be < \lat = \lat(\Om,1).  
\end{align*}
Hence $\be > \min\{  \lao (\Om, \frac{\al}{\be}), \lat(\Om,1) \}$ and using Theorem \ref{Nguyen}  with $m_p = \frac{\al}{\be}$ and $m_q=1$ we obtain that \eqref{main problem} admits a positive solution. Let $\be \le 0$. Then using Proposition \ref{PS} and Lemma \ref{MPG}, $I_+$ satisfies all the conditions of the Mountain pass theorem (see \cite[Theorem 2.1]{AmRa}). Therefore, by the Mountain pass theorem and Remark \ref{nonnegative}, \eqref{main problem} admits a nonnegative and nontrivial solution $u \in \wpso$. Further, from the strong maximum principle (Proposition \ref{SMP}), $u >0$ a.e. in $\Om$. 
\smallskip

\noi $\bm{\al < \lao,\,\be > \lat\,}$\textbf{:} Let $\al >0$. Then using $\al < \lao$ and $\be > \lat$, we get $\lat(\Om, \frac{\be}{\al})< \al < \lao(\Om,1)$. Therefore, Theorem \ref{Nguyen} with $m_p=1$ and $m_q=\frac{\be}{\al}$ yields a positive solution for \eqref{main problem}. If $\al \le 0$, then from Proposition \ref{coercive0}, we get the existence of a global minimizer $\tilde{u}$ of $I_+$, and hence using Remark \ref{nonnegative}, $\tilde{u}$ is a nonnegative solution 
 of \eqref{main problem}. Next, we show that $\tilde{u} \neq 0$ in $\Om$. Observe that, for $t>0$, $G_{\be}(t\phi_{s_2,q}) = t^q G_{\be}(\phi_{s_2,q})<0$, and using Remark \ref{Halpha}, $H_{\al}(t\phi_{s_2,q}) = t^p H_{\al}(\phi_{s_2,q})>0$. Now, if $0<t << 1$, then $I_+(t\phi_{s_2,q})<0$, which implies that $I_+(\tilde{u})<0$ and $\tilde{u} \neq 0$. Therefore, by the strong maximum principle (Proposition \ref{SMP}), $\tilde{u}>0$ a.e. in $\Om$.

\noi $\bm{\al = \lao,\,\be = \lat\,}$\textbf{:} Let \eqref{LI} violates. Then using (ii) of Proposition \ref{Nonexistence1}, we see that \eqref{main problem} admits a nonnegative solution $u \in \wpso$. Further, using the strong maximum principle (Proposition \ref{SMP}), $u>0$ a.e. in $\Om$.

\noi (ii) Suppose, there exists nonzero $c \in \R$ such that $\phi_{s_1,p}=c\phi_{s_2,q}$. We also assume that \eqref{main problem} admits a solution $u >0$ a.e. in $\Om$. Using the Picone’s inequality ((i) of Lemma \ref{picone}) and Proposition \ref{first eigenvalue}, we get
\begin{align*}
     & \iint \limits_{\rd \times \rd} |u_k(x)-u_k(y)|^{p-2}(u_k(x)-u_k(y)) \left(\frac{\phi_{s_1,p}(x)^{p}}{u_k(x)^{p-1}} - \frac{\phi_{s_1,p}(y)^p}{u_k(y)^{p-1}} \right) \, \dmuo \no \\
     & \le \iint \limits_{\rd \times \rd} \abs{\phi_{s_1,p}(x)-\phi_{s_1,p}(y)}^p \, \dmuo = \lao \intom \phi_{s_1,p}(x)^p \, \dx.
\end{align*}
Since for $x,y \in \rd$, $u_k(x)-u_k(y)=u(x)-u(y)$ the above inequality yields
\begin{equation}\label{necessary1}
\left<A_p(u), \frac{ \phi_{s_1,p}^p}{u_k^{p-1}} \right> \le \lao \intom \phi_{s_1,p}(x)^p \, \dx.
\end{equation}
We again use the Picone’s inequality ((ii) of Lemma \ref{picone}) to obtain 
\begin{equation}\label{necessary2}
    \begin{split}
         & \iint \limits_{\rd \times \rd} |u_k(x)-u_k(y)|^{q-2}(u_k(x)-u_k(y))\left(\frac{ \phi_{s_1,p}(x)^p}{u_k(x)^{p-1}}-\frac{ \phi_{s_1,p}(y)^p}{u_k(y)^{p-1}}\right) \, \dmut \\
         & \le \iint \limits_{\rd \times \rd} |\phi_{s_1,p}(x)-\phi_{s_1,p}(y)|^{q-2}(\phi_{s_1,p}(x)-\phi_{s_1,p}(y)) \left(\frac{\phi_{s_1,p}(x)^{p-q+1}}{u_k(x)^{p-q}} - \frac{\phi_{s_1,p}(y)^{p-q+1}}{u_k(y)^{p-q}} \right) \, \dmut.
    \end{split}
\end{equation}
Since $\phi_{s_1,p} \in L^{\infty}(\Om)$ ((v) of Proposition \ref{first eigenvalue}), using Lemma \ref{test functions are in space}, $u_k^{q-p}\phi_{s_1,p}^{p-q+1} \in \wpso$. Therefore, we have the following identity:
\begin{equation}\label{necessary3}
    \begin{split}
        \iint \limits_{\rd \times \rd} |\phi_{s_1,p}(x)-\phi_{s_1,p}(y)|^{q-2}(\phi_{s_1,p}(x)-\phi_{s_1,p}(y)) \left(\frac{\phi_{s_1,p}(x)^{p-q+1}}{u_k(x)^{p-q}} - \frac{\phi_{s_1,p}(y)^{p-q+1}}{u_k(y)^{p-q}} \right) \, \dmut \\ = \lat \intom \frac{\phi_{s_1,p}(x)^p}{u_k(x)^{p-q}} \, \dx.
    \end{split}
\end{equation}
Set $f_k := u_k^{q-p}\phi_{s_1,p}^p$ and $f:=u^{q-p}\phi_{s_1,p}^p$. It is easy to see that $f_k$ is increasing and $f_k \in L^1(\Om)$. Moreover, for $\ga \le p$, using $u_k(x)^{\ga-p} \ra u(x)^{\ga-p}$ a.e. in $\Om$, we get $f_k(x) \ra f(x)$  a.e. in $\Om$. Therefore, the monotone convergence theorem yields $f \in L^1(\Om)$, and $\int_{\Om} f_k(x) \, \dx \ra \int_{\Om} f(x) \, \dx$, as $k \ra \infty$.  Hence from \eqref{necessary2} and \eqref{necessary3}, we obtain
\begin{align}\label{necessary4}
     \quad  \lim_{k \ra \infty} \left<B_q(u), \frac{ \phi_{s_1,p}^p}{u_k^{p-1}} \right> \le \lat \intom \frac{\phi_{s_1,p}(x)^p}{u(x)^{p-q}} \, \dx.
\end{align}
Now since $u$ is a solution of \eqref{main problem}, taking $u_k^{1-p}\phi_{s_1,p}^{p} \in \wpso$ (by (v) of Proposition \ref{first eigenvalue}, and Lemma \ref{test functions are in space}) as a test function, 
\begin{multline}\label{necessary5}
    \left<A_p(u), \frac{ \phi_{s_1,p}^p}{u_k^{p-1}} \right> + \left<B_q(u), \frac{ \phi_{s_1,p}^p}{u_k^{p-1}} \right> = \al \intom \frac{u(x)^{p-1}}{u_k(x)^{p-1}} \phi_{s_1,p}(x)^p \, \dx + \be \intom \frac{u(x)^{q-1}}{u_k(x)^{p-1}} \phi_{s_1,p}(x)^p \, \dx. 
\end{multline}
Furthermore, for $\ga \in (1,p]$, the H\"{o}lder's inequality with the conjugate pair $(\ga,\ga')$ yields,
\begin{align*}
    \intom \frac{u(x)^{\ga-1}}{u_k(x)^{p-1}} \phi_{s_1,p}(x)^p \, \dx \le \norm{u}_{\ga}^{\ga-1} \norm{\frac{\phi_{s_1,p}^p}{u_k^{p-1}}}_{\ga} \le C(\Om, \ga) \norm{u}_{p}^{\ga-1} \norm{\frac{\phi_{s_1,p}^p}{u_k^{p-1}}}_p. 
\end{align*}
Moreover, $\frac{u^{\ga-1}}{u_k^{p-1}} \phi_{s_1,p}^p \ra u^{\ga-p} \phi_{s_1,p}^p$ a.e. in $\Om$, and the sequence $(u_k^{1-p})$ is increasing. Hence, again applying the monotone convergence theorem
\begin{align*}
    \intom \frac{u^{p-1}}{u_k^{p-1}} \phi_{s_1,p}^p \ra  \intom \phi_{s_1,p}^p \; \text{ and } \; \intom \frac{u^{q-1}}{u_k^{p-1}} \phi_{s_1,p}^p \ra \intom \frac{\phi_{s_1,p}^p}{u^{p-q}}, \text{ as } k \ra \infty.
\end{align*}
Therefore, \eqref{necessary1}, \eqref{necessary4} and \eqref{necessary5} yield
\begin{align*}
   \al \intom \phi_{s_1,p}^p + \be \intom \frac{\phi_{s_1,p}^p}{u^{p-q}} = \lim_{k \ra \infty} \left\{ \left<A_p(u), \frac{ \phi_{s_1,p}^p}{u_k^{p-1}} \right> + \left<B_q(u), \frac{ \phi_{s_1,p}^p}{u_k^{p-1}} \right> \right\} \le \lao \intom \phi_{s_1,p}^p + \lat \intom \frac{\phi_{s_1,p}^p}{u^{p-q}}.
\end{align*}
The above inequality infer that, $(\al, \be) \in ((\lao, \infty) \times (-\infty, \lat)) \cup ((-\infty, \lao) \times  (\lat,\infty)) \cup \left( \{\lao\} \times \{\lat\} \right)$. This completes our proof.\qed
\smallskip

Now we proceed to prove the existence and non-existence of positive solution for \eqref{main problem} on the line $\be = \lat$. Recall the following quantity:  
\begin{align}
    \al^*_{s_1,p} := \frac{[\phi_{s_2,q}]_{s_1,p}^p}{\norm{\phi_{s_2,q}}_p^p}.
\end{align}
Notice that, $\al^*_{s_1,p} \ge \lao$ and if \eqref{LI} holds, then $\al^*_{s_1,p} > \lao$. In the rest of this section, we assume that the condition \eqref{LI} holds.
The following lemma states that if $\al$ is smaller than  $\al^*_{s_1,p}$, then  $H_{\al}$ and $G_{\be}$ possess a different sign on $\N_{\al, \be}$.
\begin{lemma}\label{sign change}
Let $\be = \lat$ and $\al < \al^*_{s_1,p}$. Then $H_{\al}(u) < 0< G_{\be}(u)$ for every $u \in \N_{\al, \be}$.
\end{lemma}

\begin{proof}
Notice that, $G_{\be}(u) = [u]^q_{s_2,q} - \lat \norm{u^+}_q^q \ge  [u]^q_{s_2,q} - \lat \norm{u}_q^q \ge 0$ for $u \in \wpso \setminus \{0\}$. Let $u \in \N_{\al, \be}$. If $G_{\be}(u) = 0$, then we get
\begin{align*}
    \frac{[u]_{s_2,q}^q}{\norm{u}_q^q} \le \lat \le \frac{[u]_{s_2,q}^q}{\norm{u}_q^q}.
\end{align*}
By the simplicity of $\lat$ ((iv) of Proposition \ref{first eigenvalue}), $u=c\phi_{s_2,q}$ for some $c \in \R$. Hence 
\begin{align*}
    H_{\al}(u) > [u]_{s_1,p}^p - \al^*_{s_1,p} \norm{u^+}_p^p = C \left([\phi_{s_2,q}]_{s_1,p}^p -  \al^*_{s_1,p} \norm{\phi_{s_2,q}}_p^p \right)=0.
\end{align*}
On the other hand, since $u \in \N_{\al, \be}$, $H_{\al}(u) = - G_{\be}(u) = 0$, a contradiction. Therefore, we must have $G_{\be}(u) >0$. Further, since $u \in \N_{\al, \be}$, we obtain $H_{\al}(u)<0<G_{\be}(u)$.
\end{proof}

Now we are ready to prove the existence and non-existence of positive solution for $\be = \lat$.

\begin{proposition}\label{Existence2} 
For $\be = \lat$ the following hold: 
\begin{itemize}
    \item[(i)]  If $\lao < \al < \al^*_{s_1,p}$ and \eqref{LI} holds, then \eqref{main problem} admits a positive solution.
    \smallskip
    \item[(ii)] If $\al > \al^*_{s_1,p}$, then there does not exist any positive solution of \eqref{main problem}. 
\end{itemize} 
\end{proposition}

\begin{proof}
(i) We show that $d := \min \{I_+(u): u \in \N_{\al, \be}\}$ is attained. Let $(u_n)$ be the minimizing sequence in $\N_{\al,\be}$, i.e., $\big<I_+'(u_n),u_n\big>=0$ for all $n \in \mathbb{N}$ and $I_+(u_n) \ra d$ as $n \ra \infty$. From Lemma \ref{sign change},  $H_{\al}(u_n) < 0 < G_{\be}(u_n)$. 

\noi \textbf{Step 1:} This step proves the boundedness of $(u_n)$ in $\wpso$. On a contrary, suppose $[u_n]_{s_1,p} \ra \infty$, as $n \ra \infty$. Set $w_n = u_n [u_n]_{s_1,p}^{-1}$. By the reflexivity, $w_n \rightharpoonup w$ in $\wpso$ and $w_n \ra w$ in $L^p(\Om)$. Since $H_{\al}(u_n) < 0$, we have $\norm{w_n}_p^p = \norm{u_n}_p^p [u_n]^{-p}_{s_1,p} > \frac{1}{\al}$. This gives $\norm{w}_p^p \ge \frac{1}{\al}$, and hence $w \neq 0$.  Now using (i) of Remark \ref{Halpha1},
\begin{align}\label{existence1}
    \frac{p-q}{pq} G_{\be}(w_n) = \frac{I_+(u_n)}{[u_n]^q_{s_1,p}} \ra 0, \text{ as } n \ra \infty.
\end{align}
Using \eqref{existence1} we obtain 
\begin{align*}
  0 \le G_{\be}(w) \le \lowlim_{n \ra \infty} G_{\be}(w_n)=0.
\end{align*}
Therefore, $w = c\phi_{s_2,q}$ for some $c \in \R$. Further, using (i) of Remark \ref{Halpha1}, and  \eqref{existence1},
$$[\phi_{s_2,q}]_{s_1,p}^p - \al \norm{\phi_{s_2,q}}_p^p  = H_{\al}(w) \le \lowlim_{n \ra \infty} H_{\al}(w_n) = -  \lowlim_{n \ra \infty} \frac{G_{\be}(w_n)}{[u_n]_{s_1,p}^{p-q}} =0.$$ 
The above inequality yields $\al^*_{s_1,p} \le \al$, a contradiction. Therefore, $(u_n)$ must be bounded in $\wpso$. 

\noi \textbf{Step 2:} By the reflexivity, $u_n \rightharpoonup \tilde{u}$ in $\wpso$. In this step, we show $(u_n)$ converges to $\tilde{u}$ in $\wpso$. On a contrary, suppose $[u_n]_{s_1,p} \not \ra [\tilde{u}]_{s_1,p}$. If $\lim_{n \ra \infty} [u_n]_{s_1,p} < [\tilde{u}]_{s_1,p}$, then  $\lowlim_{n \ra \infty} [u_n]_{s_1,p} < [\tilde{u}]_{s_1,p}$ contradicts the weak lower semicontinuity of $[\cdot]_{s_1,p}$. Henceforth, assume that $[\tilde{u}]_{s_1,p} < \lim_{n \ra \infty} [u_n]_{s_1,p}$. Using this inequality we get $[\tilde{u}]_{s_1,p} < \lowlim_{n \ra \infty} [u_n]_{s_1,p}$ and $H_{\al}(\tilde{u}) < \lowlim_{n \ra \infty} H_{\al}(u_n) \le 0$. This implies that $\tilde{u}$ is nonzero. Now, $G_{\be}(\tilde{u}) \ge 0$, and if $G_{\be}(\tilde{u})=0$, then $\tilde{u}=c\phi_{s_2,q}$ for some $c \in \R$. Hence $H_{\al}(\phi_{s_2,q}) <0$ which implies that $\al> \al^*_{s_1,p}$, a contradiction. Therefore, $H_{\al}(\tilde{u}) < 0 < G_{\be}(\tilde{u})$. Now applying Proposition \ref{Nehari1} there exists a unique $t_{\al,\be} \in \R^+$ such that $t_{\al,\be}\tilde{u} \in \N_{\al, \be}$ and $0< I_+(t_{\al,\be} \tilde{u})= \max_{t \in \R^+} \, I_+(t\tilde{u})$. Moreover, from (ii) of Remark \ref{Halpha1}, $I_+(u_n)= \max_{t \in \R^+} \, I_+(tu_n)$. Therefore, 
\begin{align*}
    d \le I_+(t_{\al, \be} \tilde{u}) < \lowlim_{n \ra \infty} I_+(t_{\al,\be} u_n) \le  \lowlim_{n \ra \infty} I_+(u_n) = d,
\end{align*}
a contradiction. Thus, $[u_n]_{s_1,p} \ra [\tilde{u}]_{s_1,p}$ in $\R^+$. Hence from the uniform convexity of $\wpso$, $u_n \ra \tilde{u}$ in $\wpso$. 

\noi \textbf{Step 3:} In this step we prove that $\tilde{u}$ is a positive solution of \eqref{main problem}. Since $u_n \ra \tilde{u}$ in $\wpso$, we obtain $d = I_+(\tilde{u})$ and $\big<I_+'(\tilde{u}), \tilde{u} \big>=0$. Using the continuity of $H_{\al}$ and $G_{\be}$, $H_{\al}(\tilde{u}) \le 0 \le G_{\be}(\tilde{u})$. Next, we show $\tilde{u}$ is nonzero. Set $w_n = u_n [u_n]_{s_1,p}^{-1}$. Then $w_n \rightharpoonup w$ in $\wpso$. Since $H_{\al}(u_n)<0$, from the same arguments as in previous steps, $w \neq 0$ and $G_{\be}(w)>0$. Next, suppose $[u_n]_{s_1,p} \ra 0$ as $n \ra \infty$. Using $G_{\be}(w_n) \ge 0$ we get
\begin{align*}
    [w]^p_{s_1,p}-\al\norm{w}^p_p \le H_{\al}(w) \le \lowlim_{n \ra \infty} H_{\al}(w_n) = -  \lowlim_{n \ra \infty} \frac{G_{\be}(w_n)}{[u_n]_{s_1,p}^{p-q}} = -\infty.
\end{align*}
A contradiction, as $w \in L^p(\Om)$. Thus $\inf_{n \in \mathbb{N}} [u_n]_{s_1,p} >0$ and $\al \norm{\tilde{u}}_p^p \ge \lim_{n \ra \infty} [u_n]_{s_1,p}^p >0$, which implies that $\tilde{u}$ is nonzero in $\Om$, and hence $\tilde{u} \in \N_{\al, \be}$. Moreover, from Lemma \ref{sign change}, $H_{\al}(\tilde{u}) < 0 < G_{\be}(\tilde{u})$. Now, using Proposition \ref{Nehari0} and Remark \ref{nonnegative}, we conclude $\tilde{u}$ is a nonnegative solution of \eqref{main problem}. Furthermore, by Proposition \ref{SMP}, $\tilde{u}>0$ a.e. in $\Om$.

\noi (ii) Our proof uses the method of contradiction. Let $u \in \wpso$ and $u>0$ a.e. in $\Om$. From (v) of Proposition \ref{first eigenvalue} and Lemma \ref{test functions are in space}, $u_k^{q-p}\phi_{s_2,q}^{p-q+1}, u_k^{1-p}\phi_{s_2,q}^p \in \wpso$. Applying the discrete Picone’s inequality ((ii) of Lemma \ref{picone}),
\begin{equation}\label{nonexists1}
    \begin{split}
        & \iint \limits_{\rd \times \rd} |u_k(x)-u_k(y)|^{q-2}(u_k(x)-u_k(y))\left(\frac{ \phi_{s_2,q}(x)^p}{u_k(x)^{p-1}}-\frac{ \phi_{s_2,q}(y)^p}{u_k(y)^{p-1}}\right) \, \dmut \\
        & \le \iint \limits_{\rd \times \rd} |\phi_{s_2,q}(x)-\phi_{s_2,q}(y)|^{q-2}(\phi_{s_2,q}(x)-\phi_{s_2,q}(y)) \left(\frac{\phi_{s_2,q}(x)^{p-q+1}}{u_k(x)^{p-q}} - \frac{\phi_{s_2,q}(y)^{p-q+1}}{u_k(y)^{p-q}} \right) \, \dmut \\
        & = \lat \intom \frac{\phi_{s_2,q}(x)^p}{u_k(x)^{p-q}} \, \dx.
    \end{split}
\end{equation} 
The monotone convergence theorem yields $u^{q-p}\phi_{s_2,q}^p \in L^1(\Om)$ and $\int_{\Om} u_k^{q-p}\phi_{s_2,q}^p \ra \int_{\Om} u^{q-p}\phi_{s_2,q}^p$, as $k \ra \infty$. Next, we again use the Picone’s inequality ((i) of Lemma \ref{picone}), to get
\begin{equation}\label{nonexists2}
    \begin{split}
        & \iint \limits_{\rd \times \rd} |u_k(x)-u_k(y)|^{p-2}(u_k(x)-u_k(y)) \left(\frac{\phi_{s_2,q}(x)^{p}}{u_k(x)^{p-1}} - \frac{\phi_{s_2,q}(y)^p}{u_k(y)^{p-1}} \right) \, \dmuo \\
        & \le \iint \limits_{\rd \times \rd} \abs{\phi_{s_2,q}(x)-\phi_{s_2,q}(y)}^p \, \dmuo = \al^*_{s_1,p} \intom \phi_{s_2,q}(x)^p \, \dx.
    \end{split}
\end{equation}
If $u$ is a solution of \eqref{main problem}, then taking $u_k^{1-p}\phi_{s_2,q}^p$ as a test function we write
\begin{equation}\label{nonexists3}
    \begin{split}
        \left<A_p(u), \frac{\phi_{s_2,q}^{p}}{u_k^{p-1}} \right> + & \left< B_q(u), \frac{\phi_{s_2,q}^{p}}{u_k^{p-1}}  \right> \\
        & = \al \intom \frac{u(x)^{p-1}}{u_k(x)^{p-1}} \phi_{s_2,q}(x)^p \, \dx + \lat \intom \frac{u(x)^{q-1}}{u_k(x)^{p-1}} \phi_{s_2,q}(x)^p \, \dx.
    \end{split}
\end{equation}
Further, applying the monotone convergence theorem
\begin{align*}
    \intom \frac{u^{p-1}}{u_k^{p-1}} \phi_{s_2,q}^p \ra  \intom \phi_{s_2,q}^p; \quad \intom \frac{u^{q-1}}{u_k^{p-1}} \phi_{s_2,q}^p \ra \intom \frac{\phi_{s_2,q}^p}{u^{p-q}}, \text{ as } k \ra \infty.
\end{align*}
Therefore, from \eqref{nonexists1}, \eqref{nonexists2}, and \eqref{nonexists3}, we conclude 
\begin{align*}
    \al \intom \phi_{s_2,q}^p + \lat \intom \frac{\phi_{s_2,q}^p}{u^{p-q}} = \lim_{k \ra \infty} \left\{ \left<A_p(u), \frac{ \phi_{s_2,q}^p}{u_k^{p-1}} \right> + \left<B_q(u), \frac{ \phi_{s_2,q}^p}{u_k^{p-1}} \right> \right\} \le  \al^*_{s_1,p} \intom \phi_{s_2,q}^p + \lat \intom \frac{\phi_{s_2,q}^p}{u^{p-q}}.
\end{align*}
The above inequality yields $\al \le \al^*_{s_1,p}$, which is a contradiction. Thus there does not exist any positive solution for $\al > \al^*_{s_1,p}$.  
\end{proof}

\begin{remark}\label{borderline 1}
Let $\al = \al^*_{s_1,p}$ and $\be = \lat$. We assume that \eqref{LI} holds. Then observe that $I_+(u) = \frac{p-q}{pq} G_{\be}(u) \ge 0$ for every $u \in \N_{\al,\be}$, and $H_{\al}(\phi_{s_2,q})=G_{\be}(\phi_{s_2,q})=0$. Therefore, for any $t \neq 0$, we get $t \phi_{s_2,q} \in \N_{\al, \be}$ and $d=I_+(t \phi_{s_2,q}) = 0$. 
On the other hand, suppose $\phi_{s_2,q}$ is a solution of
\begin{equation}\label{P1}
    (-\Delta)_p^{s_1}u= \al^*_{s_1,p} |u|^{p-2}u \text{ in } \Om, \quad u=0 \text{ in } \rd \setminus \Om.
\end{equation}
Then $\phi_{s_2,q}$ has to change it's sign in $\Om$ (since $\al^*_{s_1,p} > \lao$), a contradiction. Thus $\phi_{s_2,q}$ does not satisfy \eqref{P1} and hence $\phi_{s_2,q}$ is not a solution of \eqref{main problem}. Thus, in this case, there does not exist any solution of \eqref{main problem} which minimizes $d$.
\end{remark}

For $\al \ge \lao$ and $\be \ge \lat$, analogously as in \cite{BoTa2} we consider the following quantity:
\begin{equation*}
    \be^{\star}(\al) = \inf \left\{ \frac{\normt^q}{\norm{u}_q^q} : u \in \wpso \setminus \{0\} \text{ and } H_{\al}(u) \le 0 \right\}.
\end{equation*}
Since $u \in \wpso \subset \wqst$, the quantity $\be^{\star}(\al) < \infty$.
\begin{proposition}
Let $\al \ge \lao$ and $\be \ge \lat$.  Assume that \eqref{LI} holds. Then $\be^{\star}(\al)$ is attained. Further, if $\al < \al^*_{s_1,p}$, then $ \be^{\star}(\al) > \lat$.
\end{proposition}

\begin{proof}
Due to the homogeneity, 
\begin{align*}
    \be^{\star}(\al) = \inf \left\{\normt^q : u \in \mathcal{M} \right\}, \text{ where } \mathcal{M} := \left\{ u \in \wpso, \norm{u}_q^q =1, \text{ and } H_{\al}(u) \le 0 \right\}.
\end{align*}
Let $(u_n)$ be a minimizing sequence for $\be^{\star}(\al)$ in $\mathcal{M}$. Suppose $[u_n]_{s_1,p} \ra \infty$. Then $H_{\al}(u_n) \le 0$ implies $\al^{\frac{1}{p}} \norm{u_n}_p \ge [u_n]_{s_1,p} \ra \infty$. Set $w_n=u_n \norm{u_n}_p^{-1}$. Then $[w_n]_{s_1,p} = [u_n]_{s_1,p} \norm{u_n}_p^{-1} \le \al^{\frac{1}{p}}$, and $w_n \rightharpoonup w$ in $\wpso$. Using $\norm{u_n}_q=1$, we get $\norm{w_n}_q \ra 0$ in $\R^+$. On the other hand, $\norm{w_n}_p=1$. Now the compact embeddings of $\wpso \hookrightarrow L^{\ga}(\Om); \ga \in [1,p]$ yields: 
\begin{align*}
  (a) \, \norm{w}_q = 0 \text{ which implies }  w=0 \text{ a.e. in } \Om; \quad (b) \, \norm{w}_p = 1.
\end{align*}
Clearly, $(a)$ and $(b)$ contradict each other. Therefore, the sequence $(u_n)$ is bounded in $\wpso$. By the reflexivity, $u_n \rightharpoonup \tilde{u}$ in $\wpso$. Further, $\tilde{u} \in \mathcal{M}$ follows from the compact embedding of $\wpso$ and weak lower semicontinuity of $H_{\al}$. Therefore, 
\begin{align*}
  \be^{\star}(\al) \le [\tilde{u}]_{s_2, q}^q \le  \lowlim_{n \ra \infty} [u_n]_{s_2, q}^q = \be^{\star}(\al).
\end{align*}
Thus $\be^{\star}(\al)$ is attained. Clearly, $\be^{\star}(\al) \ge \lat$. If $\be^{\star}(\al) = \lat$, then by the simplicity of $\lat$ ((iv) of Proposition \ref{first eigenvalue}), $\tilde{u} =c\phi_{s_2,q}$ for some $c\in \R$. Further, since $\al < \al^*_{s_1,p}$, we get $H_{\al}(\tilde{u})= CH_{\al}(\phi_{s_2,q})>0$, a contradiction to $\tilde{u} \in \mathcal{M}$. Thus $\be^{\star}(\al) > \lat$.
\end{proof}

Now we prove the existence of a positive solution for \eqref{main problem} when $\al, \be $ are larger than $\lao, \lat$ respectively. 

\begin{proposition}\label{existence 3}
Let $ \lao \le \al < \al^*_{s_1,p}$ and $\lat < \be < \be^{\star}(\al)$. Assume that \eqref{LI} holds.  Then \eqref{main problem} admits a positive solution.
\end{proposition}

\begin{proof}
We adapt the arguments as given in \cite[Theorem 2.5]{BoTa2}. As before, we will show that $d := \min \{I_+(u): u \in \N_{\al, \be}\}$ is attained. Since $\be > \lat$ and $\al <\al^*_{s_1,p}$, we have $G_{\be}(\phi_{s_2,q})<0<H_{\al}(\phi_{s_2,q})$. Then by Proposition \ref{Nehari1}, there exists a unique $t_{\al,\be} \in \R^+$ such that $0>I_+(t_{\al,\be}\phi_{s_2,q})=\min_{t \in \R^+} I_+(t\phi_{s_2,q})$. We also have $t_{\al,\be}\phi_{s_2,q} \in \N_{\al,\be}$. Therefore, $d<0$. Let $(u_n)$ be the minimizing sequence in $\N_{\al,\be}$ for $d$. Then there exists $n_0 \in \mathbb{N}$ such that $I_+(u_n)<0$ for $n \ge n_0$. Since $u_n \in \N_{\al,\be}$, using (i) of Remark \ref{Halpha1}, $G_{\be}(u_n)< 0 < H_{\al}(u_n)$ for $n \ge n_0$.

\noi \textbf{Step 1:} In this step, we show that $(u_n)$ is a bounded sequence in $\wpso$. As before, to prove this we argue by contradiction. Suppose $[u_n]_{s_1,p} \ra \infty$, as $n \ra \infty$, and set $w_n = u_n [u_n]_{s_1,p}^{-1}$. Then $w_n \rightharpoonup w$ in $\wpso$. Hence 
\begin{align*}
    H_{\al}(w) \le  \lowlim_{n \ra \infty} H_{\al}(w_n) = - \lowlim_{n \ra \infty} \frac{G_{\be}(w_n)}{[u_n]_{s_1,p}^{p-q}} = 0, \text{ and }  1- \norm{w}_p^p= \lowlim_{n \ra \infty} H_{\al}(w_n) = 0,
\end{align*}
which implies that $w \neq 0$. Therefore, from the definition,  $\be^{\star}(\al) \le [w]_{s_2,q}^q \norm{w}_q^{-q}$. Using this inequality along with $\be < \be^{\star}(\al)$, we get $G_{\be}(w) >0$. On the other hand, $G_{\be}(w) \le \lowlim_{n \ra \infty} G_{\be}(w_n) \le 0$, a contradiction. 

\noi \textbf{Step 2:} Let $u_n \rightharpoonup \tilde{u}$ in $\wpso$. This step shows that $\tilde{u}$ is a positive solution of \eqref{main problem}. First, we claim $H_{\al}(\tilde{u})>0$. On a contrary, assume $H_{\al}(\tilde{u}) \le 0$. Since $I_+(\tilde{u}) \le \lowlim_{n \ra \infty} I_+(u_n) \le d <0$, we get $\tilde{u}\neq 0.$ Hence $\be^{\star}(\al) \le [\tilde{u}]_{s_2,q}^q \norm{\tilde{u}}_q^{-q}$ and $\be < \be^{\star}(\al)$ imply $G_{\be}(\tilde{u}) > 0$. On the other hand, $G_{\be}(\tilde{u}) \le \lowlim_{n \ra \infty} G_{\be}(u_n) \le 0$, a contradiction. Therefore, $H_{\al}(\tilde{u})>0$. Further, $H_{\al}(\tilde{u})+ G_{\be}(\tilde{u}) = I_+(\tilde{u}) \le \lowlim_{n \ra \infty} I_+(u_n) \le 0$ yields $G_{\be}(\tilde{u})<0$. Now we can use Proposition \ref{Nehari1}, to get a unique $t_{\al,\be} \in \R^+$ that minimizes $I_+(t \tilde{u})$ over $\R^+$, and $t_{\al,\be} \tilde{u} \in \N_{\al,\be}$. Hence 
\begin{align*}
    d \le I_+(t_{\al,\be} \tilde{u}) = \min_{t \in \R^+} I_+(t \tilde{u}) \le I_+(\tilde{u}) \le \lowlim_{n \ra \infty} I_+(u_n) = d.
\end{align*}
Thus, $I_+(t_{\al,\be} \tilde{u})=I_+(\tilde{u})=d$ and from the uniqueness of $t_{\al, \be}$, we get $\tilde{u} \in \N_{\al,\be}$. Therefore, by Proposition \ref{Nehari0} and Remark \ref{nonnegative}, $\tilde{u}$ is a nonnegative solution of \eqref{main problem}. Further, using Proposition \ref{SMP}, $\tilde{u}>0$ a.e. in $\Om$.
\end{proof}

\begin{remark}\label{existence of solution with epsilon}
Suppose \eqref{LI} holds. We consider
\begin{align*}
    \ep_1 := \min \left\{\frac{\al^*_{s_1,p} - \lao}{2}, \frac{\be^{\star}(\al) - \lat}{2} \right\}.
\end{align*}
Then for each $\ep \in (0, \ep_1)$, using Proposition \ref{existence 3} we can conclude that (\tb{EV; $\lao+ \ep,\lat + \ep$}) admits a positive solution.
\end{remark}
Recall that, $\la^*(\theta)$ (where $\theta\in\re$) is defined as 
\begin{equation*}
    \la^*(\theta):=\sup\left\{\la\in\re:(\tb{\text{EV; $\la+\theta,\la$}})\text{ has a positive solution}\right\}.
\end{equation*}
Next, we prove some properties of the curve $\mathcal{C}:= \{(\la^*(\theta) + \theta, \la^*(\theta)): \theta \in \R\}$.  
\smallskip

\noi \textbf{Proof of Proposition \ref{threshold curve}:}
Proofs of (iii), (iv), and (vi) directly follow from \cite[Proposition~3]{BoTa1} with needful changes. So, we prove the remaining parts of the proposition.

\noi (i) Suppose $u\in\wpso$ is a positive solution of (\tb{EV; $\la+\theta,\,\la$}) for some $\la\in\re.$ For $v\in\cc(\Omega)$ with $v\geq0$, and for $k\in\mathbb{N}$, define $ \phi_{k}:= \frac{v^p}{u_{k}^{p-1}+u_{k}^{q-1}}$.
By Lemma \ref{test functions are in space}, $\phi_{k}\in\wpso.$ Using the discrete Picone's inequalities ((iii) and (iv) of Lemma \ref{picone}), we obtain
\begin{align*}
& \iint\limits_{\rd\times\rd} |u(x)-u(y)|^{p-2}(u(x)-u(y))(\phi_{k}(x)-\phi_{k}(y)) \,\dmuo
\leq\II{\rd\times\rd} |v(x)-v(y)|^p \,\dmuo, \\
& \iint\limits_{\rd\times\rd} |u(x)-u(y)|^{q-2}(u(x)-u(y))(\phi_{k}(x)-\phi_{k}(y)) \,\dmut
\leq\II{\rd\times\rd} \left|v^{\frac{p}{q}}(x)-v^{\frac{p}{q}}(y)\right|^q \,\dmut.
\end{align*}
Summing the above inequalities and by weak formulation of $u$ (where we use $\phi_{k}$ as a test function)
\begin{equation*}
    \la\intom \frac{u(x)^{p-1}+u(x)^{q-1}}{u_k(x)^{p-1}+ u_k(x)^{q-1}}v(x)^p \,\dx+\theta\intom\frac{u(x)^{p-1} v(x)^p}{u_k(x)^{p-1} + u_k(x)^{q-1}}\,\dx\leq [v]_{s_1,p}^p+\left[v^{\frac{p}{q}}\right]_{s_2,q}^q.
\end{equation*}
Further, using the monotone convergence theorem 
\begin{align*}
    \int_{\Om} \frac{u^{p-1}+u^{q-1}}{u_{k}^{p-1}+ u_{k}^{q-1}}v^p \ra \int_{\Om} v^p \, \text{ and } \, \int_{\Om} \frac{u^{p-1}v^p}{u_{k}^{p-1}+u_{k}^{q-1}} \ra \int_{\Om}  \frac{u^{p-1}v^p}{u^{p-1}+u^{q-1}}, \text{ as } k \ra \infty.
\end{align*}
This implies that 
\begin{equation}\label{part i) eq3}
    \la\intom v^p\,\dx+\min\left\{0,\theta\intom v^{p}\,\dx\right\}\leq [v]_{s_1,p}^p+\left[v^{\frac{p}{q}}\right]_{s_2,q}^q.
\end{equation}
Since $v\in\cc(\Omega)$, the R.H.S. of \eqref{part i) eq3} is a positive constant independent of $\la$ and $u$. Hence, from \eqref{part i) eq3} we conclude that $\la^*(\theta)<\infty.$

\noi (ii) \textbf{Sufficient condition:} Suppose the property \eqref{LI} holds. By remark \ref{existence of solution with epsilon}, we see that (\tb{EV; $\la^1_{s_1,p}+\ep,\,\la^1_{s_2,q}+\ep$}) admits a positive solution for $\ep>0$ small enough. From the definition of $\theta^*$, we have $(\la^1_{s_1,p}+\ep,\,\la^1_{s_2,q}+\ep)=(\la^1_{s_2,q}+\ep+\theta^*,\,\la^1_{s_2,q}+\ep)$. Hence from the definition of $\la^*(\theta^*)$,  $\la^*(\theta^*)\geq\la^1_{s_2,q}+\ep$, and $\la^*(\theta^*)+\theta^*\geq\la^1_{s_1,p}+\ep$. 

\noi \textbf{Necessary condition:} On a contrary assume that \eqref{LI} violates. This gives $\phi_{s_1,p}$ is an eigenfunction of $(-\Delta)^{s_2}_q$. Let $u$ be a positive solution of \eqref{main problem} for some $\al,\,\be\in\re.$ For $k\in\mathbb{N}$, set 
\begin{align*}
    v_k:= \frac{\phi_{s_1,p}^p}{u_k^{p-1}} \, \text{ and } \, w_k= \frac{\phi_{s_1,p}^{p-q+1}}{u_k^{p-q}}.
\end{align*}
From Lemma \ref{test functions are in space}, $v_k, w_k \in \wpso$. Using the discrete Picone's inequalities ((i) and (ii) of Lemma \ref{picone}) and Proposition \ref{first eigenvalue}, we obtain
\begin{equation}\label{part ii) Picone eq1}
    \begin{split}
        \iint \limits_{\rd\times\rd}|u(x)-u(y)|^{p-2}(u(x)-u(y))(v_k(x)-v_k(y)) \,\dmuo
        & \leq \II{\rd\times\rd} |\phi_{s_1,p}(x)-\phi_{s_1,p}(y)|^p \,\dmuo \\ 
        & = \lao \intom \phi_{s_1,p}(x)^p \, \dx,
    \end{split}
\end{equation}
and 
\begin{equation}\label{part ii) Picone eq2}
    \begin{split}
    & \iint\limits_{\rd\times\rd}  |u(x)-u(y)|^{q-2}(u(x)-u(y))(v_k(x)-v_k(y)) \,\dmut
    \\ & \leq \II{\rd\times\rd} \left|\phi_{s_1,p}(x)-\phi_{s_1,p}(y)\right|^{q-2}  \left(\phi_{s_1,p}(x)-\phi_{s_1,p}(y)\right)(w_k(x)-w_k(y)) \,\dmut 
    \\ & = \la^{1}_{s_2,q}\intom \frac{\phi_{s_1,p}(x)^p}{u_k(x)^{p-q}} \, \dx \le
    \la^{1}_{s_2,q}\intom \frac{\phi_{s_1,p}(x)^p}{u(x)^{p-q}} \, \dx,
    \end{split}
\end{equation}
where the last equality holds since $(\phi_{s_1,p}, \la^{1}_{s_2,q})$ is an eigenpair. Summing \eqref{part ii) Picone eq1}, \eqref{part ii) Picone eq2} and using $u$ is a solution of \eqref{main problem} with the test function $v_k$ we obtain
$$
\al\intom u(x)^{p-1}\frac{\phi_{s_1,p}(x)^p}{u_k(x)^{p-1}} \, \dx + \be\intom u(x)^{q-1}\frac{\phi_{s_1,p}(x)^p}{u_k(x)^{p-1}} \, \dx
\leq\la^{1}_{s_1,p}\intom\phi_{s_1,p}(x)^p \, \dx + \la^{1}_{s_2,q} \intom \frac{\phi_{s_1,p}(x)^p}{u(x)^{p-q}} \, \dx.
$$
Therefore, by the monotone convergence theorem 
$$
\al\intom\phi_{s_1,p}(x)^p\,\dx + \be \intom \frac{\phi_{s_1,p}(x)^p}{u(x)^{p-q}} \, \dx
\leq\la^{1}_{s_1,p}\intom\phi_{s_1,p}(x)^p \,\dx+\la^{1}_{s_2,q} \intom \frac{\phi_{s_1,p}(x)^p}{u(x)^{p-q}} \, \dx,
$$
a contradiction if $\al>\la^{1}_{s_1,p}$ and $\be>\la^{1}_{s_2,q}$ hold simultaneously. Thus if \eqref{LI} is violated, then there does not exist any $\be > \lat$ so that (\tb{\text{EV; $\be+\theta^*,\be$}}) admits a positive solution.

\noi (v) The proof consists of the following two cases.

\noi \textbf{Case 1:} If \eqref{LI} does not hold, then $\theta^*=\theta^*_+$ and also, $\la^*(\theta^*)\leq\la^1_{s_2,q}$ (by (ii)). Hence using the decreasing property (iv) of $\la^*(\theta)$, we get $\la^*(\theta)\leq\la^1_{s_2,q}$ for all $\theta\geq\theta^*_+.$ Therefore, the result follows in this case by using (iii).

\noi \textbf{Case 2:} Let \eqref{LI} holds. We argue by contradiction. Suppose, there exists $\theta_0\geq\theta^*_+$ such that $\la^*(\theta_0)>\la^1_{s_2,q}$. By increasing property (iv) together with (ii), we get
$\la^*(\theta_0)+\theta_0\geq\la^*(\theta^*)+\theta^*>\la^1_{s_1,p}$. By the definition of $\la^*(\theta_0)$, for any $\delta_0>0$ there exists $\delta\in[0,\delta_0)$ such that (\tb{EV; $\la^*(\theta_0)+\theta_0-\delta,\,\la^*(\theta_0)-\delta$}) admits a positive solution and we let $u$ be such solution. We choose $\de_0>0$ sufficiently small such that the following hold:
\begin{equation}\label{delta0 est}
  \la^*(\theta_0)+\theta_0-\delta_0>\la^1_{s_1,p},\text{ and }\la^*(\theta_0)-\delta_0>\la^1_{s_2,q}.  
\end{equation}
Now, by the weak formulation of $u$ and using discrete Picone's inequalities as in \eqref{part ii) Picone eq1} and \eqref{part ii) Picone eq2} (where we replace $\phi_{s_2,q}$ by $\phi_{s_1,p}$ in the test function $v_k$)
\begin{align*}
        & (\la^*(\theta_0)+\theta_0-\delta)\intom u(x)^{p-1}\frac{\phi_{s_2,q}(x)^p}{u_k(x)^{p-1}} \, \dx+(\la^*(\theta_0)-\delta)\intom u(x)^{q-1}\frac{\phi_{s_2,q}(x)^p}{u_k(x)^{p-1}}\, \dx
        \\
        & = \iint\limits_{\rd\times\rd}|u(x)-u(y)|^{p-2}(u(x)-u(y))(v_k(x)-v_k(y)) \,\dmuo \\
        & \quad + \iint\limits_{\rd\times\rd}|u(x)-u(y)|^{q-2}(u(x)-u(y))(v_k(x)-v_k(y)) \,\dmut
        \\
        & \leq\II{\rd\times\rd} |\phi_{s_2,q}(x)-\phi_{s_2,q}(y)|^p \,\dmuo +\la^{1}_{s_2,q}\intom \frac{\phi_{s_2,q}(x)^p}{u(x)^{p-q}} \, \dx.
\end{align*}
Letting $k \ra \infty$ and applying the monotone convergence theorem in the above, we obtain
\begin{equation}\label{Part v) picone eq1}
    \begin{split}
        (\la^*(\theta_0)+\theta_0-\delta)\intom \phi_{s_2,q}(x)^p \, \dx+(\la^*(\theta_0)-\delta)\intom u(x)^{q-p}\phi_{s_2,q}(x)^p \, \dx\\
\leq\II{\rd\times\rd} |\phi_{s_2,q}(x)-\phi_{s_2,q}(y)|^p \,\dmuo +\la^{1}_{s_2,q}\intom \frac{\phi_{s_2,q}(x)^p}{u(x)^{p-q}} \, \dx. 
    \end{split}
\end{equation}
Again, since $\delta<\delta_0$, we obtain from \eqref{delta0 est} that
\begin{equation}\label{part v) delta eq2}
    \begin{split}
        (\la^1_{s_2,q}+\theta_0) & \intom\phi_{s_2,q}(x)^p \,\dx+ \la^1_{s_2,q} \intom \frac{\phi_{s_2,q}(x)^p}{u(x)^{p-q}} \,\dx \\
         & < (\la^*(\theta_0)+\theta_0-\delta) \intom \phi_{s_2,q}(x)^p \dx+(\la^*(\theta_0)-\delta)\intom \frac{\phi_{s_2,q}(x)^p}{u(x)^{p-q}}\,\dx.
    \end{split}
\end{equation}
Thus, from \eqref{Part v) picone eq1} and \eqref{part v) delta eq2}
\begin{equation*}
  (\la^1_{s_2,q}+\theta_0) \intom \phi_{s_2,q}(x)^p \,\dx<\II{\rd\times\rd} |\phi_{s_2,q}(x)-\phi_{s_2,q}(y)|^p \, \dmuo, 
\end{equation*}
and this implies
$
\theta_0<\frac{[\phi_{s_2,q}]^p_{s_1,p}}{\norm{\phi_{s_2,q}}^p_p}-\la^1_{s_2,q}=\theta^*_+,
$
a contradiction to $\theta_0\geq\theta^*_+.$ This completes the proof.
\qed
\smallskip

\noi \textbf{Proof of Theorem \ref{Existence main 1}:}
(i) Suppose $\be \in (\lat, \la^*(\theta))$. From the definition of $\la^*$, there exists $\mu \in (\be, \la^*(\theta))$ such that (\tb{EV; $\mu+\theta,\mu$}) has a positive solution $\overline{u} \in \wpso$ and from Theorem \ref{bounded solution}, $\overline{u} \in L^{\infty}(\rd)$. Further, since $\mu > \be$, $\overline{u}$ is a supersolution for (\tb{EV; $\be+\theta,\be$}). Moreover, $0$ is a subsolution for (\tb{EV; $\be+\theta,\be$}). Therefore, using Proposition \ref{critical point of truncation}, $\tilde{I}$ admits a global minimizer $\tilde{u}$ in $\wpso$, and then using Proposition \ref{sub and sup} we infer that $\tilde{u} \in L^{\infty}(\rd)$ is a solution of \eqref{main problem}, satisfying $0 \le \tilde{u}(x) \le \overline{u}(x)$ a.e. in $\rd$. Next, we show that $\tilde{u}$ is nonzero. Since, $\Omega$ is a bounded open Lipschitz set then by the similar density argument given in \cite[Lemma 2.3]{BrLiPa}, there exists a sequence $\{\phi_n\} \subset \cc(\Omega)$ such that $\phi_n \ra \phi_{s_2,q}$ in $L^q(\Omega)$ and $[\phi_n]_{s_2,q} \ra [\phi_{s_2,q}]_{s_2, q}$ as $n \ra \infty$. Since $\phi_{s_2,q}>0$ in $\Omega$, without loss of generality we can choose $\phi_n$ to be non-negative in $\Omega$. As $\be>\lat$, we have  $G_{\be}(\phi_{s_2,q})<0$. Therefore, there exists $n_1 \in \mathbb{N}$ such that for $n \ge n_1$, $ [\phi_{n}]_{s_2, q}^q - \be \norm{\phi_{n}}^q_q<0.$ Further, choose $t_{n_1}>0$ small enough so that $t_{n_1} \phi_{n_1}(x) \le \overline{u}(x)$ a.e. $x \in \Omega$ and using the fact that $q<p$ we get 
\begin{align*}
    \tilde{I}(t_{n_1}\phi_{n_1}) = I(t_{n_1}\phi_{n_1})= \frac{t_{n_1}^p}{p} \left( [\phi_{n_1}]_{s_1, p}^p - \al \norm{\phi_{n_1}}^p_p \right) + \frac{t_{n_1}^q}{q} \left( [\phi_{n_1}]_{s_2, q}^q - \be \norm{\phi_{n_1}}^q_q \right)<0.
\end{align*}
Therefore, since $\tilde{u}$ is the global minimizer for $\tilde{I}$ in $\wpso$, we must have $\tilde{I}(\tilde{u})<0$, and hence $\tilde{u} \neq 0$ in $\Om$. Now, applying the strong maximum principle (Proposition \ref{SMP}), $\tilde{u}>0$ a.e. in $\Om$.

\noi (ii) Suppose $\al> \la^1_{s_1,p}$ and $\be < \la^*(\theta)$. If $\be > \lat$, then using the previous arguments, existence result holds. Now we assume $\be = \lat$. Since $\lat < \la^*(\al-\be)$ and $\la^*(\theta)$ (where $\theta = \al-\be$) is decreasing ((iv) of Proposition \ref{threshold curve}), we have $\theta < \theta^*_+$ (from (v) of Proposition \ref{threshold curve}). From the definition of $\theta^*_+$, it is easy to observe that $\theta < \theta^*_+$ is equivalent to $\al< \al^*_{s_1,p}$. Therefore, for $\be = \lat$ and $\al \in (\lao, \al^*_{s_1,p})$ using Proposition \ref{Existence2} we conclude that \eqref{main problem} admits a positive solution. 

\noi (iii) If $\be > \la^*(\al-\be)$, then from the definition of $\la^*$ we see that \eqref{main problem} does not admit any positive solution.
\qed

\noi \textbf{Proof of Theorem \ref{Existence main 2}:}
(i) Let $\theta < \theta^*_+$. From Proposition \ref{threshold curve}, $\be := \la^*(\theta) > \lat$, and  $\al := \la^*(\theta)+\theta > \lao$. From the definition of $\la^*$, there exists a sequence $(\be_n) \subset (\lat, \la^*(\theta))$, such that $\be_n \ra \be$ and (\tb{EV; $\be_n+\theta,\be_n$}) admit a sequence of positive solutions $(u_n)$ (by (i) of Theorem \ref{Existence main 1}). Now, using the similar set of arguments as given in Proposition \ref{PS}, we get  $u_n \ra \tilde{u}$ in $\wpso$. Thus, from the continuity of $I'$,  $\tilde{u}$ is a  nonnegative solution of \eqref{main problem}. 
Next, we show $\tilde{u} \neq 0$. On a contrary, assume that $\tilde{u}=0$. For each $n, k \in \mathbb{N}$, set $u_{n,k}(x)=u_n(x) + \frac{1}{k}$. From Lemma \ref{test functions are in space}, $u_{n,k}^{1-q}\phi_{s_2,q}^{q} \in \wpso$. Therefore, since $u_n$ is a solution of \eqref{main problem} 
\begin{align}\label{main 2.2}
    \left<A_p(u_n), \frac{\phi_{s_2,q}^q}{u_{n,k}^{q-1}} \right> + \left<B_q(u_n), \frac{\phi_{s_2,q}^q}{u_{n,k}^{q-1}} \right> = \al \intom u_n^{p-1} \frac{\phi_{s_2,q}^q}{u_{n,k}^{q-1}} + \be \intom u_n^{q-1} \frac{\phi_{s_2,q}^q}{u_{n,k}^{q-1}}. 
\end{align}
Using the monotone convergence theorem, $\int_{\Om} u_n^{p-1} \phi_{s_2,q}^q u_{n,k}^{1-q} \ra \int_{\Om} u_n^{p-q} \phi_{s_2,q}^q$ and $\int_{\Om} u_n^{q-1} \phi_{s_2,q}^q u_{n,k}^{1-q} \ra \int_{\Om} \phi_{s_2,q}^q$ as $k \ra \infty$. Next, from (i) of Lemma \ref{picone} and  using $u_{n,k}(x) - u_{n,k}(y) = u_n(x) - u_n(y)$, we get
\begin{equation}\label{main 2.3}
    \begin{split}
 & \left<A_p(u_n), \frac{\phi_{s_2,q}^q}{u_{n,k}^{q-1}} \right> \le \iint\limits_{\rd\times\rd} \abs{\phi_{s_2,q}(x) - \phi_{s_2,q}(y)}^q \abs{u_n(x) - u_n(y)}^{p-q} \, \dmuo. \\
    & \left<B_q(u_n), \frac{\phi_{s_2,q}^q}{u_{n,k}^{q-1}} \right> \le \iint\limits_{\rd\times\rd} \abs{\phi_{s_2,q}(x) - \phi_{s_2,q}(y)}^q \, \dmuo = \lat \norm{\phi_{s_2,q}}^q_q.
    \end{split}
\end{equation}
Further, by the H\"{o}lder inequality with the conjugate exponent $\left(\frac{p}{p-q}, \frac{p}{q}\right)$ we estimate the following inequalities:
\begin{align*}
    & \iint\limits_{\rd\times\rd} \abs{\phi_{s_2,q}(x) - \phi_{s_2,q}(y)}^q \abs{u_n(x) - u_n(y)}^{p-q} \, \dmuo \le [u_n]_{s_1,p}^{p-q} [\phi_{s_2,q}]_{s_1, p}^q,  \; \text{ and } \\
    & \int_{\Om} u_n^{p-q} \phi_{s_2,q}^q \le \norm{u_n}_p^{p-q} \norm{\phi_{s_2,q}}_p^q.
\end{align*}
Therefore, since $u_n \ra 0$ in $\wpso$, from \eqref{main 2.2} and \eqref{main 2.3} we obtain
\begin{align*}
    \be \norm{\phi_{s_2,q}}^q_q = \lim_{n \ra \infty} \lim_{k \ra \infty} \left( \left<A_p(u_n), \frac{\phi_{s_2,q}^q}{u_{n,k}^{q-1}} \right> + \left<B_q(u_n), \frac{\phi_{s_2,q}^q}{u_{n,k}^{q-1}} \right> \right) \le \lat \norm{\phi_{s_2,q}}^q_q.
\end{align*}
The above inequality yields $\be \le \lat$, a contradiction to $\be > \lat$. Therefore, $\tilde{u} \neq 0$ and from the strong maximum principle (Proposition \ref{SMP}), $\tilde{u} >0$ a.e. in $\Om$. 

\noi (ii) If $\theta > \theta^*_+$, then from (v) of Proposition \ref{threshold curve}, the problem (\tb{EV; $\la^*(\theta)+\theta,\la^*(\theta)$}) is equivalent to the problem \eqref{main problem} where $\be = \lat$ and $\al = \lat + \theta > \lat + \theta^*_+ > \al_{s_1,p}^*$. Therefore, by (ii) of Proposition \ref{Existence2}, \eqref{main problem} does not admit any positive solution.  \qed

\begin{remark}\label{borderline 2}
Let $\theta = \theta^*_+$ and \eqref{LI} holds. Then using (v) of Proposition \ref{threshold curve} and Remark \ref{borderline 1} we see that (\tb{EV; $\la^*(\theta)+\theta,\la^*(\theta)$}) does not admit any solution which minimizes $d := \min \{I_+(u): u \in \N_{\al, \be}\}$.  
\end{remark}

\begin{remark}
In this remark, we represent $\la^*(\theta)$ as a variational characterization. Let $\Omega \subset\rd$ be a bounded open set with $\C^{1,1}$ boundary $\pa \Om$. We consider the following quantity: 
\begin{align*}
    \Lambda^*(\theta) := \underset{u \in \text{int}(\C(\overline{\Om})_+)}{\sup} \underset{v \in \C(\overline{\Om})_+ \setminus \{0\}}{\inf} \frac{\left<A_p(u), v \right> + \left<B_q(u), v \right> - \theta \int_{\Om} \abs{u}^{p-2}uv}{ \int_{\Om} \left( \abs{u}^{p-2}uv + \abs{u}^{q-2}uv \right)},
\end{align*}
where $\C(\overline{\Om})_+ =  \{u \in \C(\overline{\Om}) : u \ge 0\}$ and $\text{int}(\C(\overline{\Om})_+) =  \{u \in \C(\overline{\Om})_+ : u > 0\}$. From Theorem \ref{Existence main 1}, we see that for certain ranges of $\la$, (\tb{EV; $\la+\theta,\,\la$}) admits a positive solution $u$. Further, combining Theorem \ref{bounded solution} and \cite[Corollary 2.1]{GiKuSr}, it is evident that the solution $u$ is in $\C(\overline{\Omega})$. Thus, $\text{int}(\C(\overline{\Om})_+)$ is nonempty and $\Lambda^*(\theta)$ is well defined. Now, using the same arguments as given in \cite[Proposition 5]{BoTa1} we conclude that $\la^*(\theta) = \Lambda^*(\theta)$ for every $\theta \in \R$.  
\end{remark}

\noi \textbf{Acknowledgments:}
The authors thank Prof. Vladimir Bobkov for his valuable suggestions and comments, which improved the article. N. Biswas is supported by the Department of Atomic Energy, Government of India, under project no. 12-R $\&$ D-TFR-5.01-0520. F. Sk is supported by the Alexander von Humboldt foundation.

\bibliographystyle{plainurl}

\end{document}